# The Chameleon Groups of Richard J. Thompson: Automorphisms and Dynamics

Matthew G. Brin

Department of Mathematical Sciences
State University of New York at Binghamton
Binghamton, NY 13902-6000

December 27, 1994

## Contents

Part I. Background

## 0. Introduction.

In the late 1960's Richard Thompson introduced a family of groups in connection with studies in logic. Thompson's groups have since appeared in a variety of mathematical topics: the word problem, infinite simple groups, homotopy and shape theory, group cohomology, dynamical systems and analysis. Univeral algebraic properties of these groups have been exploited in several of these connections. Thus we see that Thompson's groups have the ability to be interesting objects in many settings.

In this paper, we exploit the interaction between Thompson's groups and the subjects of ordered permutation groups and dynamical systems to analyze the automorphisms of some of the groups. Along the way we discover rigidity properties of some local structures on the real line and the circle. To show where ordered permutation groups and dynamical systems are used, and to explain the connection between the automorphisms and rigidity properties, we describe one of our results.

Thompson's groups have various representations including representations as piecewise linear bijections of the real line, the unit interval or the circle. Further, Thompson's groups can be described as all bijections that satisfy certain local properties. If we view the circle $S^1$ as $\mathbf{R}/\mathbf{Z}$ (the reals mod the integers), then one of Thompson's groups $T$ consists of all homeomorphisms $f$ from $S^1$ to itself that

(1) are piecewise linear (PL),
(2) preserve orientation,
(3) use only slopes that are integral powers of 2,
(4) have have their breaks (discontinuities of $f'$) in the set $\mathbf{Z}[\frac{1}{2}]$ (all $a/2^b$ where $a$ and $b$ are integers), and
(5) satisfy $f(\mathbf{Z}[\frac{1}{2}]) \subseteq \mathbf{Z}[\frac{1}{2}]$.

We interpret the requirement that $f$ be PL to mean that each $x$ in $S^1$ has an open neighborhood on which $x$ is (at most) the only break of $f$. From this it follows that elements of $T$ have finitely many breaks.

We prove that the outer automorphism group of $T$ is the cyclic group of order 2 generated by the homeomorphism $x \mapsto -x$ on $S^1 = \mathbf{R}/\mathbf{Z}$. The first step in this proof comes from the theory of ordered permutation groups where it is shown that any automorphism of $T$ is realized as a conjugation by some self homeomorphism of $S^1$. Thus we are reduced to the study of normalizers of $T$ in the homeomorphism group of $S^1$. The second step uses dynamical systems techniques to show that such a self homeomorphism (perhaps after composition with $x \mapsto -x$) must be in $T$. We discuss this second step further.

Properties (1)–(5) are local properties in the same way that being smooth is a local property. The group $T$ both defines and is defined by a local structure $\mathcal{G}$. There is more than one approach to defining such a structure, but the simplest is to let $\mathcal{G}$ be the set of germs of elements of $T$. A function is said to be compatible with $\mathcal{G}$ if all the germs of the function are in $\mathcal{G}$. The group $T$ then becomes the largest group of homeomorphisms that are compatible with $\mathcal{G}$. Our second step above can now be reworded to say that (up to composition with $x \mapsto -x$) a homeomorphism of $S^1$ that normalizes $T$ is compatible with $\mathcal{G}$. Since $T$ also defines $\mathcal{G}$, we can reword this further to say that (up to composition with $x \mapsto -x$) a conjugation that preserves the local structure $\mathcal{G}$ must be by an element that is compatible with $\mathcal{G}$. This last sentence states our rigidity property.

Once rigidity is established, we can apply it to analyze the automorphism groups of groups and monoids on $S^1$ other than $T$ that consist of elements that are compatible with $\mathcal{G}$. From there we can lift the local structure to the real line $\mathbf{R}$ and derive results about automorphisms of groups and monoids that are defined on $\mathbf{R}$. In particular, we analyze the automorphisms of another of Thompson's groups known as $F$ that acts either on the unit interval or $\mathbf{R}$, depending on the representation chosen.

In order to apply a single rigidity result to a large collection of groups, we develop machinery for relating the local structures that define different classes of groups. We also use this machinery to reduce the rigidity result to a much narrower statement. It is very natural to regard the elements of $T$ as generators of the structure $\mathcal{G}$. We formalize this notion and show that $\mathcal{G}$ is generated by the doubling map $x \mapsto 2x$ on $S^1$ if we introduce "flexibility" at the points in $\mathbf{Z}[\frac{1}{2}]$. Note that the doubling map is compatible with $\mathcal{G}$. With this single generator for $\mathcal{G}$ in hand, we derive rigidity from the following statement: If $h$ is an orientation preserving homeomorphism from $S^1$ to itself that preserves the set $\mathbf{Z}[\frac{1}{2}]$, and if both $h$ and $h^{-1}$ conjugate the doubling map to functions that are compatible with $\mathcal{G}$, then $h$ must be compatible with $\mathcal{G}$. Since this statement investigates conjugators of an expanding map, the use of techniques from dynamical systems is natural at this point.

We discuss another interesting aspect of Thompson's group $T$. The cyclic group of order two is also the outer automorphism group of $\mathrm{Homeo}_+(S^1)$, the full group of orientation preserving homeomorphisms of $S^1$. Thus $T$, a countable group, is imitating the behavior of a larger group. This coincidence is not isolated since both $T$ and $\mathrm{Homeo}_+(S^1)$ are simple. (The group $T$ is finitely presented, and was among the first examples of finitely presented, infinite, simple groups.) See also the remarks on Page 188 of [12] about how the cohomology of $T$ imitates that of larger groups.

Having proofs that two objects share a property makes one curious about the proofs. The more identical the proofs, the more it is likely that one object is revealing secrets about the other. (The proofs that $T$ and $\text{Homeo}_+(S^1)$ are simple can be made identical to the last step where it must be shown the commutator subgroup contains all elements that fix sets with non-empty interior. The known proofs of the last step are quite different.)

This idea has affected our choice of technique. The results that we use from ordered group theory are quite general. To apply them to Thompson's groups we only need properties (high levels of transitivity on dense subsets) that that the groups share with larger homeomorphism groups. In our use of dynamical systems techniques, we have attempted to suppress purely algebraic properties of these groups and again use only properties that they share with larger homeomorphisms groups. For example, the fact that these groups are finitely generated has been replaced in Sections 3 and 4 by the more general fact that the relevant local structures are finitely generated.

There is a rather extensive literature on Thompson's groups and their generalizations. We give some references in Section 1. However, the only properties that we need are the transitivity properties which we derive in Section 2.

We do not analyze some generalizations of Thompson's groups. For each integer $n > 1$, there is a group $F_n$ that acts on $\mathbf{R}$ of which $F_2$ is $F$. For integers $n$ and $r$ with $1 \leq r < n$, there is a group $T_{n,r}$ that acts on the circle of which $T_{2,1}$ is $T$. Among the properties that distinguish these groups from one another is the fact that elements of $F_n$ and $T_{n,r}$ all have slopes that are powers of $n$. The only automorphism groups that we analyze from these families are those of $F$ and $T$. All of our techniques apply to the groups $F_n$ and $T_{n,r}$ except the techniques in Section 11. Since we hope that the automorphism groups of all the $F_n$ and $T_{n,r}$ can be analyzed eventually, and since restricting the discussion in sections other than Section 11 would not substantially shorten the exposition (except possibly in Section 2), we have made the discussion general wherever possible.

There is also a family of groups $G_{n,r}$ that act on the Cantor set that generalize another of Thompson's groups (often referred to as $G = G_{2,1}$) that we have not mentioned above. We have no analysis of the automorphisms of any $G_{n,r}$.

**1. Statements, history and outline.**

We wish to describe certain homeomorphisms of the real line $\mathbf{R}$ and the circle $S^1$. It turns out that we will need several parametrizations of the circle, and we will use $S_r$ to denote the space $\mathbf{R}/r\mathbf{Z}$, where $\mathbf{Z}$ denotes the integers and $r$ is in $\mathbf{Z}$. This makes $S_r$ the circle of length $r$. We avoid the superscript 1 as in $S_r^1$ for purely aesthetic reasons. With this parametrization in place, we can discuss linearity, piecewise linearity and slope of functions on $S_r$ and the rationality of points in $S_r$ exactly as we would on $\mathbf{R}$.

Let $n$ be a positive integer greater than 1 and let $\mathbf{Z}[\frac{1}{n}]$ denote the set of points in $\mathbf{R}$ or $S_r$ that have the form $an^b$ where $a$ and $b$ are in $\mathbf{Z}$. Let $PL_n(\mathbf{R})$ denote the set of all homeomorphisms $f$ from $\mathbf{R}$ to itself satisfying:

(1) $f$ is piecewise linear;
(2) $f$ is orientation preserving;
(3) all slopes of $f$ are integral powers of $n$;
(4) the "breaks" of $f$ (discontinuities of $f'$) are in a discrete subset of $\mathbf{Z}[\frac{1}{n}]$; and
(5) $f(\mathbf{Z}[\frac{1}{n}]) \subseteq \mathbf{Z}[\frac{1}{n}]$.

It is easy to verify that every $f \in PL_n(\mathbf{R})$ satisfies $f(\mathbf{Z}[\frac{1}{n}]) = \mathbf{Z}[\frac{1}{n}]$ and that $PL_n(\mathbf{R})$ is closed under composition and inversion and is thus a group. The term "break" goes back to [8] where the simplicity of certain groups of PL hoemeomorphisms of $\mathbf{R}$ is proven.

We are interested in various subgroups of $PL_n(\mathbf{R})$ and groups on the circles $S_r$ that are related to $PL_n(\mathbf{R})$. Any unsupported statements that we make about the groups that we define are not needed in this paper. Justifications can be found in the references that we give later in this section. Given a homeomorphism $h$ from a topological space $X$ to itself, we define the support of $f$ to be the set of points $x \in X$ for which $f(x) \neq x$. We define $BPL_n(\mathbf{R})$ to be the elements of $PL_n(\mathbf{R})$ whose support is a bounded subset of $\mathbf{R}$, and we refer to the elements of $BPL_n(\mathbf{R})$ as the elements of $PL_n(\mathbf{R})$ with bounded support.

We let $F_n$ be those elements $f$ of $PL_n(\mathbf{R})$ that are translations by multiples of $(n-1)$ near $\pm\infty$ in that there are integers $i$ and $j$ and a real $M$ so that $f(x) = x + i(n-1)$ for all $x > |M|$ and $f(x) = x + j(n-1)$ for all $x < -|M|$. We have $BPL_n(\mathbf{R}) \subseteq F_n$. The technicalities in the definition of $F_n$ come from our choice of representation. (A fact not needed in this paper is that $F_n$ is isomorphic to the elements of $PL_n(\mathbf{R})$ with support in the unit interval $[0,1]$. The isomorphism can be achieved as a conjugation by a suitably chosen PL homeomorphism (with infinitely many breaks) from the open unit interval to the real line.) Note that (5) follows from

(5′) $f(x) \in \mathbf{Z}[\frac{1}{n}]$ for some $x \in \mathbf{Z}[\frac{1}{n}]$

so that (5) follows from (1)–(4) for $F_n$ and $BPL_n(\mathbf{R})$.

We let $T_{n,r}$ be those homeomorphisms from $S_r$ to itself that satisfy (1)–(5) above. (We will not need the fact that $T_{n,r}$ and $T_{n,r'}$ are isomorphic if $r$ and $r'$ are congruent mod $n-1$.) In imitation of $BPL_n(\mathbf{R})$, we let $BT_{n,r}$ be the subgroup of $T_{n,r}$ that is generated by those elements of $T_{n,r}$ whose support misses an open interval of $S_r$. That is, $BT_{n,r}$ is generated by elements of $T_{n,r}$ that are the identity on some open interval of $S_r$. It will be seen later that $BT_{n,r}$ is sometimes a proper subgroup of $T_{n,r}$.

We also consider elementary extensions of these groups. If we replace (2) with

(2′) $f$ is orientation preserving or orientation reversing,

then we get groups $\widetilde{PL_n}(\mathbf{R})$ and $\widetilde{T}_{n,r}$ that respectively contain $PL_n(\mathbf{R})$ and $T_{n,r}$ as subgroups of index 2. We also consider certain monoids. We let $\overline{PL}_n(\mathbf{R})$ be the set of functions $f$ whose domain and range are open sets of $\mathbf{R}$ that satisfy (1)–(5) where (2) is to be interpreted locally and (4) is to require that breaks of $f$ occur in $\mathbf{Z}[\frac{1}{n}]$ and be discrete in the domain of $f$. The same definition for functions from open sets of $S_r$ to open sets in $S_r$ gives $\overline{T}_{n,r}$. These are monoids if we make composition defined on the largest domain possible and if we include the empty function as a zero element. Note that the degree $n$ map $x \mapsto nx$ is an element of $\overline{T}_{n,r}$. Lastly we consider $\widetilde{\overline{PL}}_n(\mathbf{R})$ and $\widetilde{\overline{T}}_{n,r}$ which are the obvious two-fold extensions of $\overline{PL}_n(\mathbf{R})$ and $\overline{T}_{n,r}$ obtained by also allowing functions that are locally orientation reversing. We do not allow functions that preserve orientation on some open sets and reverse orientation on others.

At this point we have defined the largest and smallest objects that we will look at. Our results will be be about all groups and certain monoids $G$ for which $BPL_n(\mathbf{R}) \subseteq G \subseteq \widetilde{\overline{PL}}_n(\mathbf{R})$ or $BT_{n,r} \subseteq G \subseteq \widetilde{\overline{T}}_{n,r}$. Our "smallest" groups are $BPL_n(\mathbf{R})$ and $BT_{n,r}$ since, as we will see in Section 2, they have all the transitivity properties that we will need.

We will need a certain amount of control on a monoid if we are to say anything about its automorphisms. Note that the invertible elements of a monoid $G$ (those $a \in G$ for which there is a $b \in G$ with $ab$ and $ba$ the identity in $G$) form a group. This group will have particular importance since our results for monoids will be consequences of our results for groups. This is emphasized in the next definitions. If $X$ is a topological space, then a local homeomorphism $f$ on $X$ is a function from a union of opens sets $U_\alpha$ of $X$ into $X$ so that the restriction of $f$ to each $U_\alpha$ is an embedding onto an open set in $X$. It is not required that $f$ be injective. If $G$ is a monoid of local homeomorphisms on $X$ and $f$ is a local homeomorphism on $X$, then we say that $f$ is compatible with the invertible elements of $G$ if for every $x$ in the domain of $f$ there is an invertible $a \in G$ that agrees with $f$ on an open set about $x$. We say that $G$ is determined by its invertible elements if $G$ consists exactly of all local homeomorphisms of $X$ that are compatible with the invertible elements of $G$. Note that the invertible elements of $G$ are automatically compatible with the invertible elements of $G$. Note also that while the restrictions of invertible elements of $G$ to open subsets of $X$ will all be in $G$, these will not be all the elements of $G$. As examples, $\overline{PL}_n(\mathbf{R})$ and $\overline{T}_{n,r}$ are determined by their invertible elements.

We are interested in automorphisms of the objects that we have defined. We get certain automorphisms by conjugation. Let $G$ be a group or monoid for which $BPL_n(\mathbf{R}) \subseteq G \subseteq \widetilde{\overline{PL}}_n(\mathbf{R})$. We let Homeo($\mathbf{R}$) to be the set of all self homeomorphisms of $\mathbf{R}$ and define $N(G)$, the normalizer of $G$ in Homeo($\mathbf{R}$), to be the set of all $h \in$ Homeo($\mathbf{R}$) for which $hGh^{-1} = G$. For $G$ with $BT_{n,r} \subseteq G \subseteq \widetilde{\overline{T}}_{n,r}$, we define $N(G)$ similarly with respect to Homeo($S_r$). We have a homomorphism $\Phi : N(G) \to \operatorname{Aut}(G)$ where for an $h \in N(G)$, we have $(\Phi(h))f = hfh^{-1}$.

We can now give our main result. Recall that $F_2 = F$ and $T_{2,1} = T$.

THEOREM 1. *Let $G$ be a group or monoid for which $BPL_n(\mathbf{R}) \subseteq G \subseteq \widetilde{PL_n}(\mathbf{R})$ or $BT_{n,r} \subseteq G \subseteq \widetilde{T}_{n,r}$. If $G$ is not a group, then assume $G$ is determined by its invertible elements. Then*

(i) $\Phi : N(G) \to Aut(G)$ *is an isomorphism.*

*For $n = 2$, we further have*

(ii) $N(G) \subseteq \widetilde{PL_2}(\mathbf{R})$ *or $N(G) \subseteq \widetilde{T}_{2,1}$ whichever applies,*
(iii) *the containment in (ii) is equality if $G$ is one of $BPL_2(\mathbf{R})$, $PL_2(\mathbf{R})$, $\widetilde{PL}_2(\mathbf{R})$, $\overline{PL}_2(\mathbf{R})$, $\widetilde{\overline{PL}}_2(\mathbf{R})$, $BT_{2,1}$, $T_{2,1}$, $\widetilde{T}_{2,1}$, $\overline{T}_{2,1}$, or $\widetilde{\overline{T}}_{2,1}$, and*
(iv) *if $N_o(F_2)$ represents the index 2 subgroup of $N(F_2)$ of orientation preserving elements, then there is a short exact sequence*

$$1 \to F_2 \to N_o(F_2) \to T_{2,1} \times T_{2,1} \to 1.$$

REMARK: The appearance of $T_{2,1} \times T_{2,1}$ in (iv) above is easy to explain. Elements of $F_2$ are translations by integers near $\pm\infty$. Thus near each of $-\infty$ and $+\infty$, elements of $F_2$ commute with functions that are lifts of homeomorphism of $S_1$. Those homeomorphisms of $S_1$ that lift to $PL_2(\mathbf{R})$ are elements of $T_{2,1}$. For an $n$ other than 2, the statement of (iv) would be more complicated. See Section 6.

We can give our "intermediate" result. It is actually not weak enough to be truly intermediate. Let $\nu_2(x) = 2x$. If we regard this as a function on $S_1$ then it is an element of $\overline{T}_{2,1}$.

THEOREM 2. *Let $h : S_1 \to S_1$ be an orientation preserving homeomorphism for which $h(\mathbf{Z}[\frac{1}{2}]) = \mathbf{Z}[\frac{1}{2}]$, and assume that $h\nu_2 h^{-1}$ and $h^{-1}\nu_2 h$ are both in $\overline{T}_{2,1}$. Then $h \in T_{2,1}$.*

The requirement $h(\mathbf{Z}[\frac{1}{2}]) = \mathbf{Z}[\frac{1}{2}]$ is an annoying technical hypothesis. We can give a cosmetically improved statement. It is "equivalent" to Theorem 2. Note that the unique fixed point of $\nu_2$ on $S_1$ is the element 0 in $S_1$.

COROLLARY 2.1. *Let $h : S_1 \to S_1$ be an orientation preserving homeomorphism that fixes 0 in $S_1$, and assume that $h\nu_2 h^{-1}$ and $h^{-1}\nu_2 h$ are both in $\overline{T}_{2,1}$. Then $h \in T_{2,1}$.*

The true intermediate result is the following.

THEOREM 3. *Let $h : S_1 \to S_1$ be an orientation preserving homeomorphism for which $h(\mathbf{Z}[\frac{1}{2}]) = \mathbf{Z}[\frac{1}{2}]$, and assume that $h\nu_2 h^{-1}$ and $h^{-1}\nu_2 h$ are both in $\overline{T}_{2,1}$. Then $h$ is PL.*

The proofs of Theorems 1 and 2 from Theorem 3 will be given in Part III as well as the "equivalence" of Corollary 2.1 and Theorem 2. It follows from the arguments in Part III that if Theorem 3 is known for $\nu_n$, $S_{n-1}$, $\mathbf{Z}[\frac{1}{n}]$ and $T_{n,n-1}$ for values of $n$ other than 2, then Theorems 1 and 2 will hold for these objects and these values of $n$ as well. The proof of Theorem 3 occupies Part IV.

A separate statement of a rigidity property is deferred to Part III. See Theorem 6.5. This is because the statement needs the terminology of Section 3, and because Part III is a convenient place where all the ingredients needed in the statement can be pulled together.

We give a brief outline of the paper.

There have been many facts discovered about the groups that we consider. However the only properties that we need are the transitivity properties. These are presented in Section 2. Thus the current paper can be read without reference to the (rather extensive) literature on the groups considered. The transitivity properties are crucial to all other sections of the paper.

Provision (i) of Theorem 1 follows directly from a result in the theory of ordered permutation groups due to McCleary and Rubin [18] and a short lemma that we prove that covers the monoid case. Thus we are in a position of needing only to verify hypotheses. This is all done in Section 5. Once (i) of Theorem 1 is proven, we are reduced to the study of normalizers rather than automorphisms.

We consider many groups of homeomorphisms defined on various domains. One of our tasks is to reduce the scope of the investigation by showing that all the normalizers that we wish to study are locally like certain conjugators of a very small set of maps. We build some general machinery for this in Section 3 that applies to groups of homeomorphisms that are defined by local properties. Section 4 applies this to the Thompson groups where we identify the small set of maps as the maps $\nu_n$ on the circles $S_{n-1}$.

In Part III, we show that all PL normalziers can be analyzed. From the results of Section 4, this gives the implication that all normalizers (and thus automorphisms) are analyzed if a certain set of conjugators of the maps $\nu_n$ contain only PL homeomorphisms. This is done in the argument that Theorem 3 implies Theorems 1 and 2.

Part IV contains the proof of Theorem 3. Here it must be shown that any homeomorphism $h$ from $S_1$ to itself that has $h\nu_2 h^{-1}$ and $h^{-1}\nu_2 h$ in $\overline{T}_{2,1}$ must be PL. Since $\nu_2$ is an expanding map, we extract a large amount of information about $h$ from our knowledge of $h\nu_2 h^{-1}$. We use this information to do a fairly direct calculation of $h^{-1}\nu_2 h$. It is then possible to tie together enough properties $h^{-1}\nu_2 h$ with properties of $h$ to show that $h$ must be PL under the hypotheses of the theorem. The difficulty in the argument is the amount of information that must be recorded before the calculation of $h^{-1}\nu_2 h$ can proceed. The bulk of the effort in Part IV goes into the development of structures that keep track of the data.

Part V contains examples that justify hypotheses and phases of the arguments. This section comes at the end of the paper since the examples are most easily described after the machinery in Part IV is developed.

Some unpublished work by others on the automorphisms of Thompson's groups has been available to the author. All have made the observation in the proof of Statement 6.2 that the main problem is to prove that normalizers are PL. Also, all have made the observation not needed or mentioned in this paper that a normalizer that is PL on some open interval must be PL. A preprint version [11] of [12] that considered the question of automorphisms gave us the idea used in Part IV of studying Markov partitions for the conjugates of the maps $\nu_n$. The author is also indebted to Dennis Sullivan for suggesting that it might be useful to look at the non-linearity of the conjugates of the maps $\nu_n$.

The most extensive work on the material considered in this paper is the unpublished notes [1] of Bieri and Strebel. (A very small sampling of material from [1] is given in an appendix to [22].) The overlap between [1] and the current paper outside of Part IV is large. Section 2 of the current paper gives a special case of the analysis in [1], and Sections 3, 4 and 6 generalize and repackage ideas in [1]. The theorem of [18] that we quote in Section 5 is a generalization of a theorem in [1] which applies to $\mathbf{R}$ which in turn is a generalization of a theorem in [17] which applies to more transitive actions. Most of the examples in Section 13 have their equivalents in [1] but the machinery we use to build the examples is different and the properties are therefore easier to verify. Items from Part IV found in [1] include a formula that gives $h$ from $h\nu_2 h^{-1}$ and an invariant based on summing break values (which are defined below in Section 9). Missing from [1] is most of the analysis of Part IV.

We state two results from [1] that complement the results of the current paper. Considered in [1] are groups $G(I; A, P)$ of all orientation preserving, PL self homeomorphisms of an interval $I$ in $\mathbf{R}$ with slopes in a multiplicative subgroup $P$ of the positive reals and breaks in a *finite* subset of an additive $\mathbf{Z}P$ module $A$ in $\mathbf{R}$ (with action coming from multiplication in $\mathbf{R}$). Also considered is the subgroup $B(I; A, P)$ of those elements of $G(I; A, P)$ with support in some compact subset of the interior of $I$. It is proven in [1] that if $G$ is a group with $B(I; A, P) \subseteq G \subseteq G(I; A, P)$ where $P$ is not cyclic, then all automorphisms of $G$ are realized as conjugations by PL self homeomorphisms of $I$. This result uses the density of such a $P$ in the positive reals. It is also proven in [1] that automorphisms of $G(I; A, P)$ are realized by PL self homeomorphisms of $I$ if $I$ is unbounded and by PL self homeomorphisms of $I$ with finitely many breaks if $I = \mathbf{R}$. This result uses the fact that restriction to an unbounded end gives a surjection onto a group of affine functions. The structures of the automorphism groups are also examined in [1].

The result, Theorem 5.1, that we use from [18] realizes automorphisms as conjugations. An earlier version of this result is found in [17], and is used there to analyze the automorphism groups of other homeomorphism groups. It is shown that the full group of PL homeomorphisms on $\mathbf{R}$ has trivial outer automorphism group as does the group of PL homeomorphisms on $\mathbf{R}$ with finitely many breaks and also the largest group of homeomorphisms of $\mathbf{R}$ all of whose elements are (not necessarily continuously) differentiable. See [17] for more details and older references. A special case of Theorem 5.1 for Thompson's group $T$ appears in [12]. The proof in [12] makes use of algebraic properties of $T$ including its simplicity.

The literature on Thompson's groups is very scattered. Because of this, many letters from different alphabets have been used to refer to Thompson's groups. Thompson's groups also admit a wide variety of representations, so the difficulty of recognition extends beyond the multiplicity of notation. We will not give a complete history, we will not sort out the notation, and we will not list all the representations. We will try to give enought references to enable the reader to get to the literature on the various topics in which these groups have appeared.

A paper that relates some of the representations, that gives a key to some of the notation used, and that gives a brief history is [3]. This is also done in [7].

For relations with the word problem and for the early history of Thompson's groups see [19].

The connection to infinite simple groups is discussed in [3] and more recently in [21]. See also [22].

Thompson's group $F$ was rediscovered (by Dydak and Minc and also independently by Freyd and Heller) in connection with questions in homotopy/shape theory. This is discussed in [10] and [6]. It is in this topic that the universal algebraic properties of the groups are derived and exploited. The final paper [16] on the homotopy question is the first to investigate the homological properties of one of Thompson's groups.

The paper [6] is the start of a systematic investigation of the (co)homology properties of Thompson's groups and their generalizations. More recent papers are [3], [4], [5], [12], [22] and [9]. The (co)homology of structures resembling those discussed in Section 3 are investigated in [13]. An acyclic extension of the infinite braid group by $BPL_2(\mathbf{R})$ is constructed in [15].

The paper [12] considers the group $T$ from a dynamic point of view.

In analysis, the existence of Thompson's group $F$ demonstrates that either there is a non-amenable finitely presented group with no free subgroup on two generators or there is a finitely presented amenable group that is not elementary amenable. This is discussed in [7] where it is shown that $F$ is not elementary amenable and in [2] and [7] where it is shown that $F$ contains no free subgroup on two generators. The monograph [7] also contains a version of Thompson's original proof that $T$ is simple and a proof of Thurston's observations that $T$ and $F$ are isomorphic to groups of piecewise projective homeomophisms on the circle and the interval respectively. The projective aspects of $T$ are also discussed in [14].

## 2. Transitivity properties of Thompson's groups.

This is the only section in which a restriction to $n = 2$ would result in substantial simplification.

We will concentrate on those properties of Thompson's groups that are needed for our analysis. It turns out that the only properties we will need are transitivity and related properties (ability to approximate), and not various algebraic properties such as being finitely generated or simple.

Since translations have slope 1 and no breaks, all translations by elements of $\mathbf{Z}[\frac{1}{n}]$ are in $PL_n(\mathbf{R})$. Similarly, all rotations by elements of $\mathbf{Z}[\frac{1}{n}]$ are in $T_{n,r}$. Thus $PL_n(\mathbf{R})$ and $T_{n,r}$ are transitive on $\mathbf{Z}[\frac{1}{n}]$. This is false in general for $BPL_n(\mathbf{R})$, $F_n$ and $BT_{n,r}$. We will look at this and at higher levels of transitivity that these groups exhibit on subsets of $\mathbf{R}$ and $S_r$. The main point is that these groups are very "flexible." Any element of $\text{Homeo}_+(\mathbf{R})$ and $\text{Homeo}_+(S_r)$ can be approximated locally by an element of $BPL_n(\mathbf{R})$ or $BT_{n,r}$.

The main obstruction to strong transitivity of all the groups that we consider is the following which corresponds to "casting out nines" when $n = 10$. It is also the reason that considering arbitrary $n$ in this section is more involved than just considering $n = 2$.

LEMMA 2.1. *There is a surjective ring homomorphism $\phi_n : \mathbf{Z}[\frac{1}{n}] \to \mathbf{Z}_{n-1}$.*

PROOF: Let

$$(2\text{-}1) \qquad x = \sum_{i=1}^{j} a_i n^{b_i}$$

where each $a_i$ and $b_i$ are in $\mathbf{Z}$. Let $k = \sum a_i$ and define $\phi_n(x) = k \mod (n-1)$. We first show that this is well defined. We can replace any one $a_i n^{b_i}$ by $(na_i)n^{b_i-1}$ and $k$ is altered by $na_i - n = (n-1)a_i$. With a finite number of alterations of this type, any two sums (2-1) that give $x$ can be reduced to two expressions $x = k_1 n^{b_1}$ and $x = k_2 n^{b_2}$. Another flurry of the same kind of steps will make the two expressions identical. This shows that $\phi_n(x)$ is well defined. Now any two elements of $\mathbf{Z}[\frac{1}{n}]$ can be represented as multiples of the same integral power of $n$. This makes it clear that $\phi_n$ is a ring homomorphism.

LEMMA 2.2. *There is a homomorphism $d_n : PL_n(\mathbf{R}) \to \mathbf{Z}_{n-1}$ so that for all $f \in PL_n(\mathbf{R})$ and $x \in \mathbf{Z}[\frac{1}{n}]$, $\phi(f(x)) = \phi(x) + d_n(f)$.*

PROOF: We must show that $d_n(f) = \phi(f(x)) - \phi(x)$ is independent of $x \in \mathbf{Z}[\frac{1}{n}]$. Equivalently, we must show for all $x$ and $y$ in $\mathbf{Z}[\frac{1}{n}]$, the truth of

$$\phi(f(x)) - \phi(x) = \phi(f(y)) - \phi(y),$$

or equivalently

(2-2) $$\phi(f(x)) - \phi(f(y)) = \phi(x) - \phi(y).$$

We can find points $x = x_0 < x_1 < \cdots < x_k = y$ so that every break of $f$ in $[x, y]$ is one of the $x_i$. For $0 \le i < k$, there are $b_i \in \mathbf{Z}$ and $c_i \in \mathbf{Z}[\frac{1}{n}]$ so that $x \in [x_i, x_{i+1}]$ implies $f(x) = n^{b_i}x + c_i$. Since $\phi$ is a ring homomorphism and $\phi(n)$ is the multiplicative identity in $\mathbf{Z}_{n-1}$, we get

$$
\begin{aligned}
(2\text{-}3) \quad \phi(f(x_{i+1})) - \phi(f(x_i)) &= \phi(f(x_{i+1}) - f(x_i)) \\
&= \phi(n^{b_i}x_{i+1} - n^{b_i}x_i) \\
&= \phi(n^{b_i})\phi(x_{i+1} - x_i) \\
&= \phi(x_{i+1} - x_i) \\
(2\text{-}4) \quad &= \phi(x_{i+1}) - \phi(x_i)
\end{aligned}
$$

and (2-2) follows by summing (2-3) and (2-4). To see that $d_n$ is a homomorphism, note that $x \in \mathbf{Z}[\frac{1}{n}]$ implies $g(x) \in \mathbf{Z}[\frac{1}{n}]$ and

$$
\begin{aligned}
d_n(fg) &= \phi(f(g(x))) - \phi(x) \\
&= [\phi(f(g(x))) - \phi(g(x))] + [\phi(g(x)) - \phi(x)] \\
&= d_n(f) + d_n(g).
\end{aligned}
$$

Let $\Delta_n$ be the kernel of $\phi_n : \mathbf{Z}[\frac{1}{n}] \to \mathbf{Z}_{n-1}$. We have $\Delta_n = \mathbf{Z}[\frac{1}{n}]$ if and only if $n = 2$. We now know that an $f \in PL_n(\mathbf{R})$ induces a specific rotation on the cosets of $\Delta_n$ in $\mathbf{Z}[\frac{1}{n}]$. The elements of the kernel of $d_n : PL_n(\mathbf{R}) \to \mathbf{Z}_{n-1}$ are those elements that induce the trivial permutation on the cosets of $\Delta_n$. The oddities in the definition of $F_n$ are now partly explained by the next corollary. Recall that $BPL_n(\mathbf{R}) \subseteq F_n$.

COROLLARY 2.2.1. *We have $F_n \subseteq \ker(d_n)$ and thus elements of $BPL_n(\mathbf{R})$ and $F_n$ preserve the cosets of $\Delta_n$ in $\mathbf{Z}[\frac{1}{n}]$ as sets.*

PROOF: Near $\pm\infty$ the action of elements of $F_n$ is to add multiples of $(n-1)$. This preserves cosets of $\Delta_n$.

We now see how transitive $BPL_n(\mathbf{R})$ is on $\Delta_n$.

LEMMA 2.3. *Let $x_1 < x_2 < \cdots < x_n$ and $y_1 < y_2 < \cdots < y_k$ be elements of $\Delta_n$. Then there is an $f \in BPL_n(\mathbf{R})$ so that $f(x_i) = y_i$ for all $i$ with $0 \leq i \leq k$.*

REMARK: We refer to the property in the lemma by saying that $BPL_n(\mathbf{R})$ acts order $k$-transitively or o-$k$-transitively on $\Delta_n$ for all $k$. Any $G$ with $BPL_n(\mathbf{R}) \subseteq G \subseteq \widetilde{PL}_n(\mathbf{R})$, is o-$k$-transitive on $\Delta_n$ for all $k$ as well.

PROOF: By adding an extra point from $\Delta_n$ to each sequence sufficiently far below $x_1$ and $y_1$ and another from $\Delta_n$ sufficiently far above $x_k$ and $y_k$, we can assume that $x_1 = y_1$ and that $x_k = y_k$. If for each $i$ with $0 \leq i < k$, we find a homeomorphism from $[x_i, x_{i+1}]$ to $[y_i, y_{i+1}]$ that satisfies the local properties required of functions in $PL_n(\mathbf{R})$, then we can piece these functions together with the identity on $(-\infty, x_1] \cup [x_k, \infty)$ to get the required function. Thus it suffices to consider four points $a < b$, $c < d$ in $\Delta_n$ and build a homeomorphism from $[a, b]$ to $[c, d]$ that satisfies the properties of functions in $PL_n(\mathbf{R})$.

Since all four points are in $\Delta_n \subseteq \mathbf{Z}[\frac{1}{n}]$, there is an integral power of $n$ that evenly divides $b-a$ and $d-c$. Thus we can divide the intervals $[a, b]$ and $[c, d]$ into subintervals of length an integral power of $n$. We know that the number of these intervals in each of $[a, b]$ and $[c, d]$ is 0 modulo $(n-1)$. If the number of intervals is not the same for $[a, b]$ and $[c, d]$, then we can increase the number of intervals for one of them as much as we want by repeatedly subdividing random subintervals into $n$ equally spaced smaller intervals. Each time this is done, the number of intervals increases by $(n-1)$. This preserves the property that the lengths of all the intervals are a (perhaps varying) integral power of $n$. Once the number of subintervals in $[a, b]$ and $[c, d]$ are the same, then a function can be built to map each interval in the subdivision of $[a, b]$ affinely to the corresponding interval in the subdivision of $[c, d]$.

The action is not $k$-transitive on all of $\mathbf{R}$. However, the set $\Delta_n$ is dense in $\mathbf{R}$, so we get close. If $I$ and $J$ are two closed intervals in $\mathbf{R}$, then we write $I < J$ if every point in $I$ is less than every point in $J$. We say that a group acting on $\mathbf{R}$ is acting approximately o-$k$-transtively on $\mathbf{R}$ if whenever points $x_1 < x_2 < \cdots < x_k$ and closed intervals $I_1 < I_2 < \cdots < I_k$ with non-empty interiors are given in $\mathbf{R}$ then there is an element of the group that takes each $x_i$ into $I_i$.

COROLLARY 2.3.1. *Let $G$ be a group with $BPL_n(\mathbf{R}) \subseteq G \subseteq \widetilde{PL}_n(\mathbf{R})$. Then $G$ acts approximately o-$k$-transitively on $\mathbf{R}$ for all $k$.*

The situation is slightly different on the circle. There is an ambiguity of reprentatives for elements of $\mathbf{Z}[\frac{1}{n}]$ in $S_r = \mathbf{R}/r\mathbf{Z}$. If $r \not\equiv 0 \mod (n-1)$ then there is no homomorphism from $\mathbf{Z}[\frac{1}{n}]$ in $S_r$ to $\mathbf{Z}_{n-1}$ but rather one to $\mathbf{Z}_a$ where $a$ is the greatest common divisor of $r$ and $(n-1)$. (In fact it turns out that $T_{n,r}$ is simple when $a = 1$. See [3].) There is also the problem of not having a linear order on $S_r$. This is solved by referring to the counterclockwise orientation on $S_r$. Since functions in $T_{n,r}$ lift to functions in $PL_n(\mathbf{R})$ and since functions in $PL_n(\mathbf{R})$ can be built in pieces, we leave it as an easy exercise for the reader to verify the following.

LEMMA 2.4. *Let $x_1, x_2, \ldots, x_k$ be different points in $S_r$ that are arranged in counterclockwise order on $S_r$ in the order listed, and let $I_1, I_2, \ldots, I_k$ be pairwise disjoint closed intervals with non-empty interiors in $S_r$ that are arranged in counterclockwise order on $S_r$ in the order listed. Then there is an element in $BT_{n,r}$ that takes each $x_i$ into $I_i$.*

We will need the following approximation lemma.

LEMMA 2.5. *Given $f \in \ker(d_n)$ and $-\infty < a < b < \infty$, there is an element $g$ in $BPL_n(\mathbf{R})$ that aggrees with $f$ on $[a, b]$.*

PROOF: By decreasing $a$ and increasing $b$ we can assume that $a$ and $b$ are in $\Delta_n$. By the hypothesis on $f$ we know that $f(a)$ and $f(b)$ are in $\Delta_n$. By Lemma 2.3, there is a $g_1$ in $BPL_n(\mathbf{R})$ that carries $a$ to $f(a)$ and $b$ to $f(b)$. Now we let $g$ agree with $f$ on $[a, b]$ and with $g_1$ elsewhere.

LEMMA 2.6. *Let $\Gamma$ be generated by translations $x \mapsto x + a$ with $a \in \Delta_n \cap (-\frac{1}{2}, \frac{1}{2})$ and elements of $BPL_n(\mathbf{R})$ with support in $(0, 1)$. Then $\Gamma$ contains $BPL_n(\mathbf{R})$.*

PROOF: The set of translations $x \mapsto x + a$ with $a \in \Delta_n$ is a topological group and the set of translations in the statement of the lemma is a neighborhood of the identity. Thus $\Gamma$ contains all translations by elements of $\Delta_n$. This puts all elements of $BPL_n(\mathbf{R})$ with support in an interval of length 1 in $\Gamma$. It is now an easy matter to generate all elements of $BPL_n(\mathbf{R})$. For example, let $f$ have support in an interval $J$ and let $b$ be an element of $\Delta_n$ very near the midpoint of $J$. We can assume that $J$ is not contained in $I = (b - \frac{1}{2}, b + \frac{1}{2})$. We can find functions $g_1, g_2, \ldots, g_k$ whose supports are in intervals of length less than 1 that are contained in the support of $f$, and whose composition $g = g_1 g_2 \ldots g_k$ takes $b$ to $f(b)$. Now $g^{-1} f_1$ has its support in $J$ and a fixed point at $b$. There is an $h$ with support in a small neighborhood of $b$ so that $hg^{-1}f_1$ is fixed near $b$ and has support contained in the union of two disjoint closed intervals with lengths less that half the length of $J$. The rest is clear.

Let $p_r : \mathbf{R} \to \mathbf{R}/r\mathbf{Z}$ be projection. If $r \equiv 0 \mod (n-1)$, then $\Delta_n = p_r^{-1}(p_r(\Delta_n))$. If not, then we get a proper containment $\Delta_n \subsetneq p_r^{-1}(p_r(\Delta_n))$. The next lemma shows that the group $BT_{n,r}$ depends on the image of $\Delta_n$ under $p_r$.

LEMMA 2.7. *The group $BT_{n,r}$ consists of those elements of $T_{n,r}$ that preserve $p_r(\Delta_n)$ and it is generated by the rotations $x \mapsto x + a$ with $a \in p_r(\Delta_n) \cap (-\frac{1}{2}, \frac{1}{2})$ and elements of $BT_{n,r}$ with support in $(0, 1)$.*

PROOF: We first argue that $BT_{n,r}$ contains the rotations by elements of $p_r(\Delta_n)$. Given a rotation $\rho(x) = x + a$ with $a \in \Delta_n \cap (-\frac{1}{2}, \frac{1}{2})$ (we bring in $p_r$ later), Lemma 2.5 builds a generator $g$ of $BT_{n,r}$ that imitates a portion of $\rho$ so that $g^{-1}\rho$ is another generator of $BT_{n,r}$. Now as in the previous lemma, $BT_{n,r}$ contains all rotations by all elements of $\Delta_n$. However, on $S_r$, these are just rotations by elements of $p_r(\Delta_n)$. Any element in $T_{n,r}$ that preserves $p_r(\Delta_n)$ can be composed with a rotation by an element of $p_r(\Delta_n)$ so that the composition has a fixed point. Now a composition with a generator of $BT_{n,r}$ gives a fixed interval and a new generator of $BT_{n,r}$. The rest of the proof repeats arguments used here and in the proof of the previous lemma.

Note that the previous lemma shows that $BT_{n,r} = T_{n,r}$ when $p_r(\Delta_n) = \mathbf{Z}[\frac{1}{n}]$ in $S_r$. This will happen when $r$ is relatively prime to $(n-1)$. As mentioned above, $T_{n,r}$ is simple when $r$ is relatively prime to $(n-1)$. We will not need these facts. We will need to know that $p_{n-1}(\Delta_n) = \Delta_n$ in $S_{n-1}$. This will make $BT_{n,n-1}$ a particularly important group.

PART II. GENERALITIES

## 3. Germs, half germs, germ functions and germ generators.

The restrictions (1)–(5) on $PL_n(\mathbf{R})$ and related groups are local in nature. In this section, we make this statement more precise and derive some consequences. On the way, we develop some generalities about locally defined classes of functions. We start with the particular and move to the general.

We start with $PL_n(\mathbf{R})$ as typical of the groups that we study. Consider an $f \in PL_n(\mathbf{R})$. From (4) and (5), we know that each point in $\mathbf{R}$ lies in a closed interval $J$ with endpoints in $\mathbf{Z}[\frac{1}{n}]$ so that $f|_J$ is the restriction to $J$ of an affine function with slope $n^k$, $k \in \mathbf{Z}$, and $y$-intercept $b \in \mathbf{Z}[\frac{1}{n}]$. Specifically, $f|_J(x) = n^k x + b$.

The behavior of $x \mapsto n^k + b$ can be better understood if we express elements of $\mathbf{R}$ as expansions in base $n$. That is, each $x \in \mathbf{R}$ is representable as

$$(3\text{-}1) \qquad x = \sum_{i=-\infty}^{\infty} a_i n^i, \ \ a_i \in \{0, 1, 2, \ldots, n-1\}$$

where $a_i = 0$ for all sufficiently large $i$. The expression is unique for elements of $\mathbf{R} - \mathbf{Z}[\frac{1}{n}]$. For elements of $\mathbf{Z}[\frac{1}{n}]$, there are two possible representations, one with all $a_i = 0$ for all sufficiently small $i$, and one with all $a_i = n - 1$ for all sufficiently small $i$. We adopt the convention that all elements of $\mathbf{Z}[\frac{1}{n}]$ are represented so that all $a_i = 0$ for sufficiently small $i$. Note that this makes the sum of the $a_i$ definable for elements of $\mathbf{Z}[\frac{1}{n}]$. We can use this to identify the elements of $\Delta_n$ as those elements of $\mathbf{Z}[\frac{1}{n}]$ for which the sum of the $a_i$ is congruent to 0 modulo $(n-1)$.

We use the usual notation

(3-2) $$x = a_j a_{j-1} \ldots a_0 . a_{-1} a_{-2} \ldots$$

to symbolize the expression in (3-1) where all $a_i = 0$ for $i > j$. We will refer to the $a_i$ as letters (or numerals) in the alphabet $\{0, 1, \ldots, n-1\}$ and we will refer to the expressions (3-2) as words that are infinite to the right with a single position between two letters marked by "." the $n$-ary point.

The act of taking $x$ to $n^k x$ shifts the word (3-2) representing $x$ to the left $k$ places with respect to the $n$-ary point (with net shift to the right if $k$ is negative). Since $b \in \mathbf{Z}[\frac{1}{n}]$, the act of adding $b$ is to alter a finite number of letters in (3-2) at the left end of the word. From this we see that if $x = \ldots a_1 a_0 . a_{-1} a_{-2} \ldots$ and $y = \ldots b_1 b_0 . b_{-1} b_{-2} \ldots$ are in the same orbit under $PL_n(\mathbf{R})$, then

(3-3) $$\text{there are integers } j \text{ and } k \text{ so that for all } i < j, \, a_{i-k} = b_i.$$

Conversely, if (3-3) is satisfied, then $y - n^k x$ is in $\mathbf{Z}[\frac{1}{n}]$ and $(x, y)$ is on the graph of $t \mapsto n^k t + (y - n^k x)$. We can now define a relation on $\mathbf{R}$ by saying that $x \sim y$ with $x$ and $y$ as represented above if (3-3) is satisfied. This is an equivalence relation whose equivalence classes are the orbits of $PL_n(\mathbf{R})$. Note that we have shown that the orbits of $PL_n(\mathbf{R})$ are exactly the orbits of

$$\text{Aff}_n(\mathbf{R}) = \{x \mapsto n^k x + b \mid k \in \mathbf{Z}, \, b \in \mathbf{Z}[\tfrac{1}{n}]\}.$$

We can now discuss germs. We start with the "smaller" group $\text{Aff}_n(\mathbf{R})$. Let $x \sim y$ in $\mathbf{R}$. The set of germs $\text{Aff}_n(\mathbf{R})_{x,y}$ will be the set of equivalence classes of those $f \in \text{Aff}_n(\mathbf{R})$ that have $f(x) = y$ under the relation $f \sim g$ if there is an open set about $x$ on which $f$ and $g$ agree. Since an affine map is determined by its behavior on an open set, we see that $\text{Aff}_n(\mathbf{R})_{x,y}$ is simply the set of affine maps in $\text{Aff}_n(\mathbf{R})$ that take $x$ to $y$. These in turn are indexed by the set of integers $k$ for which $t \mapsto n^k t + (y - n^k x)$ is in $\text{Aff}_n(\mathbf{R})$. That is, $y - n^k x \in \mathbf{Z}[\frac{1}{n}]$. If there is more than one such $k$, say $k_1$ and $k_2$, then $x - n^{k_1 - k_2} x$ is in $\mathbf{Z}[\frac{1}{n}]$, $x$ is in the rationals $\mathbf{Q}$ and the expression $x = \ldots a_1 a_0 . a_{-1} a_{-2} \ldots$ ends in a periodic pattern with period dividing $k_1 - k_2$. Conversely, if $x$ is in $\mathbf{Q}$, then the expression $x = \ldots a_1 a_0 . a_{-1} a_{-2} \ldots$ ends in a periodic pattern of minimum period $d$ and the slopes that appear in $\mathrm{Aff}_n(\mathbf{R})_{x,y}$ are exactly those slopes of the form $n^{k+di}$ for all $i \in \mathbf{Z}$, where $k$ is some integer determined modulo $d$ by $x$ and $y$. From this, we have that $\mathrm{Aff}_n(\mathbf{R})_{x,y}$ has one element if $x$ is irrational and that the slopes $\mathrm{Aff}_n(\mathbf{R})_{x,y}$ are all integral powers of $n$ if $x \in \mathbf{Z}[\frac{1}{n}]$. (Note that there are rationals outside of $\mathbf{Z}[\frac{1}{n}]$ with this last property if $n > 2$.)

If we now turn to $PL_n(\mathbf{R})$, we get sets of germs $PL_n(\mathbf{R})_{x,y}$ that are identical to those of $\mathrm{Aff}_n(\mathbf{R})_{x,y}$ except where $x \in \mathbf{Z}[\frac{1}{n}]$. For $x \in \mathbf{Z}[\frac{1}{n}]$, behavior to the right of $x$ and to the left of $x$ are independent and we get all integral powers of $n$ as slopes to the right and all integral powers of $n$ as slopes to the left.

We can now generalize. Let $X$ be a topological space. Let $\mathcal{H}(X)$ be the set of homeomorphisms from open sets in $X$ to open sets in $X$. For each $x$ and $y$ in $X$, let $\mathcal{H}(X)_{x,y}$ be the set of elements of $\mathcal{H}(X)$ that take $x$ to $y$. Let $\mathcal{G}(X)_{x,y}$ be the set of germs at $x$ of the elements of $\mathcal{H}(X)_{x,y}$. That is $\mathcal{G}(X)_{x,y}$ is the set of equivalence classes in $\mathcal{H}(X)_{x,y}$ under the relation that declares two elements related if they agree on some neighborhood of $x$.

We can compose germs in that elements of $\mathcal{G}(X)_{x,y}$ compose with those of $\mathcal{G}(X)_{y,z}$ to give elements of $\mathcal{G}(X)_{x,z}$. Also, each $\mathcal{G}(X)_{x,x}$ contains the germ of the identity. If we regard $\mathcal{G}(X)$ as a function defined on $X \times X$, then we have that $\mathcal{G}(X)$ is a category whose objects are the elements of $X$ and whose morphism sets are the sets $\mathcal{G}(X)_{x,y}$. Each element in $\mathcal{G}(X)_{x,y}$ has an inverse in $\mathcal{G}(X)_{y,x}$ so that the two possible compositions give the identities of $\mathcal{G}(X)_{x,x}$ and $\mathcal{G}(X)_{y,y}$. Thus all morphisms are isomorphisms and $\mathcal{G}(X)$ is a groupoid. We will deal with substructures of this groupoid. We will have no need to consider topological groupoids.

Let $\mathcal{G}$ be a function on $X \times X$ written so that $\mathcal{G}_{x,y}$ represents the value of $\mathcal{G}$ at $(x,y)$ and assume that for each $(x,y) \in X \times X$ we have $\mathcal{G}_{x,y} \subseteq \mathcal{G}(X)_{x,y}$. If we further assume that $\mathcal{G}$ is a groupoid, then we call $\mathcal{G}$ a groupoid of invertible germs on $X$. We will usually omit the word invertible. It is standard to arrive at a groupoid of germs from a pseudogroup (here, the set of $\mathcal{G}$-compatible maps defined below). We will not bother to do this, although we are not really avoiding that much work. See for example [20].

We say that $f : U \to X$, where $U$ is an open subset of $X$, is $\mathcal{G}$-compatible at $x \in U$ if there is a $g$ representing a class in $\mathcal{G}_{x,f(x)}$ that agrees with $f$ on some neighborhood of $x$. We say that $f$ is $\mathcal{G}$-compatible if it is $\mathcal{G}$-compatible at $x$ for each $x \in U$. We let $F(\mathcal{G})$ be the set of all $\mathcal{G}$-compatible functions, we let $H(\mathcal{G})$ be the set of all $\mathcal{G}$-compatible homeomorphisms from $X$ to $X$, and we let $O(\mathcal{G})$ be the orbits of $F(\mathcal{G})$. We have that $O(\mathcal{G})$ is just the set of equivalence classes in which two points $x$ and $y$ are related when $\mathcal{G}_{x,y}$ is not empty. Note that every element of $F(\mathcal{G})$ is a local homeomorphism in that every point in its domain has an open neighborhood that is carried homeomorphically onto an open set in $X$.

Note that we have not made enough assumptions to guarantee that $H(\mathcal{G})$ is not empty. This can be fixed by requiring that each $\mathcal{G}_{x,x}$ be non-empty. This will mean that $\mathcal{G}_{x,x}$ contains at least the germ of the identity function and this will make the identity function on $X$ an element of $F(\mathcal{G})$ and thus $H(\mathcal{G})$. Equivalently we could require that the identity function on $X$ be in $F(\mathcal{G})$. When $H(\mathcal{G})$ contains the identity function on $X$, we will say that the groupoid $\mathcal{G}$ contains the identity. Note that this is implied if we require that for each $x \in X$, there is a $y$ so that $\mathcal{G}_{x,y}$ is not empty.

To imitate the relationship between $PL_n(\mathbf{R})$ and $\operatorname{Aff}_n(\mathbf{R})$, we need to split germs in half. We now work on the real line $\mathbf{R}$ and we let $\mathcal{G}$ be a groupoid of germs on $\mathbf{R}$. Given a set of germs $\mathcal{G}_{x,y}$ and a function $f$ taking $x$ to $y$, we say that $f$ is $\mathcal{G}^+$-compatible at $x$ if there is an $\epsilon > 0$ so that $f$ and some representative in $\mathcal{G}_{x,y}$ agree on $[x, x+\epsilon)$. Similarly, $f$ is $\mathcal{G}^-$-compatible if $f$ and a representative in $\mathcal{G}_{x,y}$ agree on some $(x-\epsilon, x]$. If $A$ is a union of orbits in $O(\mathcal{G})$, then we say that a local homeomorphism $f: U \to \mathbf{R}$ is $(\mathcal{G}; A)$-compatible if $f$ is $\mathcal{G}$-compatible for all $x \in U - A$ and if $f$ is both $\mathcal{G}^+$-compatible and $\mathcal{G}^-$-compatible for each $x \in A$. Note that we must insist in advance that $f$ be a local homeomorphism since without making that requirement explicit, an increasing germ from $\mathcal{G}^+_{x,f(x)}$ might be combined with a decreasing germ from $\mathcal{G}^-_{x,f(x)}$ or the reverse. We denote the groupoid of germs of all $(\mathcal{G}; A)$-compatible functions by $(\mathcal{G}; A)$. We let $F(\mathcal{G}; A)$ be the set of all $(\mathcal{G}; A)$-compatible maps, and we let $H(\mathcal{G}; A)$ be the set of all $(\mathcal{G}; A)$-compatible homeomorphisms. The constructions in this paragraph can be compared to Paragraph 1.6 in [13].

Note that $(x, y) \mapsto \operatorname{Aff}_n(\mathbf{R})_{x,y}$ defines a groupoid of germs on $\mathbf{R}$. To avoid ambiguities of notation, we use $\mathcal{A}_n(\mathbf{R})$ rather than $\operatorname{Aff}_n(\mathbf{R})$ to denote this groupoid. Now if $A = \mathbf{Z}[\frac{1}{n}]$, then

$$H(\mathcal{A}_n(\mathbf{R})) = \operatorname{Aff}_n(\mathbf{R}),$$
$$H(\mathcal{A}_n(\mathbf{R}); A) = PL_n(\mathbf{R}), \text{ and}$$
$$F(\mathcal{A}_n(\mathbf{R}); A) = \overline{PL}_n(\mathbf{R}).$$

We now consider the circle $S_r$ and its covering map, the projection $\mathbf{R} \to \mathbf{R}/r\mathbf{Z}$. If $\mathcal{G}$ is a groupoid of germs on $S_r$, then we can lift $\mathcal{G}$ to a groupoid of germs on $\mathbf{R}$. One way to define this is to let $\mathcal{P}$ be the set of germs of the covering projection $\mathbf{R} \to \mathbf{R}/r\mathbf{Z}$ and let $\mathcal{P}^{-1}$ be the inverses of the elements of $\mathcal{P}$. Then the lift of $\mathcal{G}$ to $\mathbf{R}$ is the set of all allowable compositions in $\mathcal{P}^{-1}\mathcal{G}\mathcal{P}$. We can build a groupoid of germs $\mathcal{A}_n(S_r)$ on $S_r$ that lifts to $\mathcal{A}_n(\mathbf{R})$ on $\mathbf{R}$. We can be precise about this too. Here we are going from $\mathbf{R}$ to $S_r$, so let $\mathcal{G}$ be a groupoid of germs on $\mathbf{R}$, and let $s: \mathbf{R} \to \mathbf{R}$ be defined by $s(x) = x + r$. We say that $\mathcal{G}$ is $s$-invariant if the germ of $s^q\alpha s^{-p}$ at $x+pr$ is in $\mathcal{G}_{x+pr,y+qr}$ for every germ $\alpha$ in $\mathcal{G}_{x,y}$ and every $p$ and $q$ in $\mathbf{Z}$. If $\mathcal{G}$ contains the identity, then the $s$-invariance of $\mathcal{G}$ is equivalent to the requirement that $s$ be in $H(\mathcal{G})$. This is seen by noting that the germ of $s$ at $x$ is the germ of $s^1\alpha s^0$ at $x$ where $\alpha$ is the germ of the identity at $x$. Assuming that $\mathcal{G}$ is $s$-invariant, we can now let the projection of $\mathcal{G}$ onto $S_r$ be the set of all allowable compositions in $\mathcal{P}\mathcal{G}\mathcal{P}^{-1}$. Since $s$ is in $\mathcal{A}_n(\mathbf{R})$, we let $\mathcal{A}_n(S_r)$ be the projection of $\mathcal{A}_n(\mathbf{R})$ onto $S_r$. The set of $\mathcal{A}_n(S_r)$-compatible maps that are defined on all of $S_r$ is the monoid generated by the map $\nu_n(x) = nx$ and the rotations by elements of $\mathbf{Z}[\frac{1}{n}]$. If $A = \mathbf{Z}[\frac{1}{n}]$, then
$$H(\mathcal{A}_n(S_r); A) = T_{n,r}, \text{ and}$$
$$F(\mathcal{A}_n(S_r); A) = \overline{T}_{n,r}.$$

We now consider how little information is needed to specify a set of germs. Consider the map $\nu_n(x) = nx$ on the circle $S_1$. We represent the elements of $S_1$ as words (3-2) where all $a_i = 0$ for positive $i$. The action of $\nu_n$ on these words is to shift the entire word one letter to the left and to discard (replace with zero) the letter that passes to the left of the $n$-ary point. From this we see that if $x$ and $y$ satisfy (3-3), then there are non-negative integers $a$ and $b$ so that $\nu_n^a(x) = \nu_n^b(y)$. The converse is clear and we see that (3-3) captures not the orbits of $\nu_n$ but the symmetric, transitive closure of the orbits of $\nu_n$. The non-negative powers of $\nu_n$ form a semigroup (the word monoid leads to a phonetically unacceptable result) and so the germs of $\nu_n$ could be called a semigroupoid. Since $\nu_n$ is a local homeomorphism, we have that the germs of $\nu_n$ are all invertible and in fact the groupoid of germs $\mathcal{A}_n(S_1)$ is the smallest groupoid that contains the germs of $\nu_n$. (We do not need this exact statement, and the suspicious reader can wait for the more precise statement and proof of Lemma 6.4.) We refer to this by saying that $\nu_n$ generates $\mathcal{A}_n(S_1)$. We can abuse notation a bit and regard $\nu_n(x) = nx$ as a function on $\mathbf{R}$. If we let $\tau(x) = x+1$ on $\mathbf{R}$, then we have that $\mathcal{A}_n(\mathbf{R})$ is generated by the pair of functions $\nu_n$ and $\tau$. (This will not be needed.)

We become more formal. Let $X$ be a topological space. Let $\mathcal{F}$ be a collection of local homeomorphisms defined on open subsets of $X$. Note that the set of germs of $\mathcal{F}$ is a subset of the collection of germs in $\mathcal{G}(X)$. Thus compositions of germs of $\mathcal{F}$ and inverses of germs of $\mathcal{F}$ exist in $\mathcal{G}(X)$. Standard arguments show that the collection of finite (allowable) compositions of germs from $\mathcal{F}$ and their inverses is a groupoid and is thus the smallest groupoid that contains the germs of $\mathcal{F}$. We refer to this groupoid as the groupoid of germs generated by $\mathcal{F}$. The orbits of the groupoid of germs generated by $\mathcal{F}$ are the equivalence classes of the symmetric, transitive closure of the relation that is defined by declaring $x \sim y$ if there is a non-negative integer $a$ and an $f \in \mathcal{F}$ so that $y = f^a(x)$. The following is clear.

LEMMA 3.1. *Let $X$ be a topological space and let $\mathcal{F}$ and $\mathcal{F}_1$ be collections of local homeomorphisms defined on open subsets of $X$.*

  (i) *If $\mathcal{F} \subseteq \mathcal{F}_1$, then the groupoid of germs generated by $\mathcal{F}$ is contained in the groupoid of germs generated by $\mathcal{F}_1$.*
  (ii) *If $\mathcal{F}_c$ is the set of finite compositions of elements of $\mathcal{F}$, or if $\mathcal{F}_c$ is a group and $\mathcal{F}$ is a set of generators of $\mathcal{F}_c$, then the groupoid of germs generated by $\mathcal{F}_c$ equals the groupoid of germs generated by $\mathcal{F}$.*

We need a definition to state the next lemma. If $\mathcal{G}$ is a groupoid of germs on a topological space $X$, then we say that a homeomorphism $h : X \to X$ preserves $\mathcal{G}$ under conjugation if $hfh^{-1}$ is in $F(\mathcal{G})$ for all $f \in F(\mathcal{G})$. Note that this implies that $hfh^{-1}$ is in $H(\mathcal{G})$ for all $f \in H(\mathcal{G})$.

LEMMA 3.2 (CONSISTENCY). *Let $\mathcal{F}$ be a family of local homeomorphisms on a topological space $X$ and let $\mathcal{G}$ be the groupoid of germs generated by $\mathcal{F}$. Let $h : X \to X$ be a homeomorphism and assume that $h$ conjugates each element of $\mathcal{F}$ into $F(\mathcal{G})$. Then $h$ preserves $\mathcal{G}$ under conjugation. If in addition, $X$ is one of $S_r$ or $\mathbf{R}$, $A$ is a union of orbits in $O(\mathcal{G})$, $h(A) \subseteq A$, and $h$ conjugates each element of $\mathcal{F}$ into $F(\mathcal{G}; A)$, then $h$ preserves $(\mathcal{G}; A)$ under conjugation.*

PROOF: We give the argument that applies when $X$ is $\mathbf{R}$ or $S_r$, and let the reader omit the irrelevant parts for the simpler case. Let $f$ be in $F(\mathcal{G}; A)$. If $x \notin A$, then $h^{-1}(x) \notin A$ and the germ of $f$ at $h^{-1}(x)$ is a composition of a finite number of germs or their inverses from elements of $\mathcal{F}$. If $x \in A$, then $h^{-1}(x)$ might be in $A$ or not. If not, then the previous statement about $f$ applies. If $h^{-1}(x) \in A$, then the two half germs of $f$ at $h^{-1}(x)$ are compositions of corresponding half germs from $\mathcal{F}$. (Which half germ of each factor goes into each half germ of $f$ depends on the order in which the factors are orientation reversing or preserving.) The conjugate of the composition is the composition of the conjugates. Since the conjugate of each factor is known to have its germ or two half germs in $(\mathcal{G}; A)$, we have that $hfh^{-1}$ is in $F(\mathcal{G}; A)$.

The notation $(\mathcal{G}; A)$ refers not only to a groupoid, but also emphasizes its construction. We could use a simpler notation (a single letter, for example) to refer to the groupoid so constructed. For the remainder of the section, we will not be concerned with the construction of groupoids and notations like $(\mathcal{G}; A)$ will not be needed.

The hypothesis of the next corollary is implied by a stronger condition that we will introduce later.

COROLLARY 3.2.1. *Let $\mathcal{G}$ be a groupoid of germs on a topological space $X$ and assume that $H(\mathcal{G})$ generates $\mathcal{G}$. If $h : X \to X$ is a homeomorphism for which $hH(\mathcal{G})h^{-1} \subseteq H(\mathcal{G})$, then $h$ preserves $\mathcal{G}$ under conjugation.*

If $\mathcal{G}$ is a groupoid of germs, then $\mathcal{G}_{x,x}$ is the set of germs at $x$ of those elements of $F(\mathcal{G})$ that have $x$ as a fixed point. The following standard observation will be needed in several places.

LEMMA 3.3. *Let $X$ be a topological space and let $\mathcal{G}$ be a groupoid of germs on $X$ that contains the identity. Then each $\mathcal{G}_{x,x}$ is a group. If $h : X \to X$ is a homeomorphism and $h$ and $h^{-1}$ preserve $\mathcal{G}$ under conjugation, then the function from $\mathcal{G}_{x,x}$ to $\mathcal{G}_{h(x),h(x)}$ taking each $\alpha \in \mathcal{G}_{x,x}$ to the germ of $h\alpha h^{-1}$ at $h(x)$ is an isomorphism.*

We can illustrate this with examples that we have been discussing. With $\mathcal{A}_n(\mathbf{R})$ the groupoid of germs generated by $\mathrm{Aff}_n(\mathbf{R})$ and $A = \mathbf{Z}[\frac{1}{n}]$, then $(\mathcal{A}_n(\mathbf{R}); A)_{x,x}$ is isomorphic to $\mathbf{Z} \times \mathbf{Z}$ if $x \in \mathbf{Z}[\frac{1}{n}]$, to $\mathbf{Z}$ if $x \in \mathbf{Q} - \mathbf{Z}[\frac{1}{n}]$, and to $\{1\}$ if $x \notin \mathbf{Q}$. Thus any homeomorphism $h : \mathbf{R} \to \mathbf{R}$ that normalizes $F(\mathcal{A}_n(\mathbf{R}); A)$ must preserve the three sets, $\mathbf{Z}[\frac{1}{n}]$, $\mathbf{Q} - \mathbf{Z}[\frac{1}{n}]$ and $\mathbf{R} - \mathbf{Q}$.

Let $X$ be a topological space, let $\mathcal{G}$ be a groupoid of germs on $X$, and let $H(X)$ be the set of self homeomorphisms of $X$. We say that $g \in H(X)$ is obtained from $f \in H(X)$ by $\mathcal{G}$-rearrangement if for every $x \in X$, the germ of $g$ at $x$ is that of $\alpha f \beta$ where $\alpha$ and $\beta$ are germs in $\mathcal{G}$ that vary with $x$. Note that this relation does not seem to be symmetric. This relation becomes more symmetric when we consider normalizers. We let $N(\mathcal{G})$ be the set of $h \in H(X)$ so that both $h$ and $h^{-1}$ preserve $\mathcal{G}$ under conjugation. Note that $H(\mathcal{G})$ is a normal subgroup of $N(\mathcal{G})$.

LEMMA 3.4 (REARRANGEMENT). *Let $X$ be a topological space and let $\mathcal{G}$ be a groupoid of germs on $X$ that contains the identity. If $f$ and $g$ in $H(X)$ are in the same left coset of $H(\mathcal{G})$ in $H(X)$ or the same right coset of $H(\mathcal{G})$ in $H(X)$, then each is a $\mathcal{G}$-rearrangment of the other. If $h$ is in $N(\mathcal{G})$, then any $\mathcal{G}$-rearrangement of $h$ is in $N(\mathcal{G})$ and elements of $N(\mathcal{G})$ are in the same coset of $H(\mathcal{G})$ in $N(\mathcal{G})$ if and only if they are $\mathcal{G}$-rearrangements of each other.*

PROOF: If $g = fH(\mathcal{G})$, then $g = fh$ with $h \in H(\mathcal{G})$. This and a similar observation for right cosets finishes the first claim. If $h$ is in $N(\mathcal{G})$, $g$ is a $\mathcal{G}$-rearrangement of $h$ and $f$ is in $F(\mathcal{G})$, then the composition of germs

$$gfg^{-1} = (\alpha h\beta)f(\alpha' h\beta')^{-1}$$

shows that $gfg^{-1}$ is in $F(\mathcal{G})$ as well. We know that elements of the same coset of $H(\mathcal{G})$ in $N(\mathcal{G})$ are $\mathcal{G}$-rearrangements of each other. Now if $g$ and $h$ in $N(\mathcal{G})$ satisfy $g = \alpha h\beta$ at each $x$ for some germs $\alpha$ and $\beta$ that vary with $x$, then $gh^{-1} = \alpha h\beta h^{-1}$ at each $x$, and $h \in N(\mathcal{G})$ implies that $gh^{-1} \in F(\mathcal{G})$ and thus in $H(\mathcal{G})$, and $g$ and $h$ are in the same coset of $H(\mathcal{G})$ in $N(\mathcal{G})$.

We borrow a property from the study of foliations. Let $\mathcal{G}$ be a groupoid of germs on $\mathbf{R}$. We say that $\mathcal{G}$ is interpolating if for every pair of germs $\alpha \in \mathcal{G}_{a,b}$ and $\beta \in \mathcal{G}_{c,d}$ with both $\alpha$ and $\beta$ increasing, with $a < c$ and with $b < d$, there is an $f \in F(\mathcal{G})$ with $\alpha$ the germ of $f$ at $a$ and $\beta$ the germ of $f$ at $c$. This is condition $(\gamma)$ in [20] and is called *germ-connected* in [13].

LEMMA 3.5. *Let $\mathcal{G}$ be an interpolating groupoid of germs on $\mathbf{R}$ that contains the identity. Then for every pair of germs $\alpha \in \mathcal{G}_{a,b}$ and $\beta \in \mathcal{G}_{c,d}$ with both $\alpha$ and $\beta$ increasing, with $a < c$ and with $b < d$, there is an $h \in H(\mathcal{G})$ with $\alpha$ the germ of $h$ at $a$ and $\beta$ the germ of $h$ at $c$.*

PROOF: There is an $f \in F(\mathcal{G})$ with $\alpha$ the germ of $f$ at $a$ and $\beta$ the germ of $f$ at $c$. For some $x$ with $x < a$ and $x < b$ there is a $g_1 \in F(\mathcal{G})$ with $\alpha$ the germ of $g_1$ at $a$ and with the germ of $g_1$ at $x$ the germ of the identity. For some $y$ with $y > c$ and $y > d$ there is a $g_2 \in F(\mathcal{G})$ with $\beta$ the germ of $g_2$ at $c$ and with the germ of $g_2$ at $y$ the germ of the identity. Now let

$$h(t) = \begin{cases} t & t \leq x, \\ g_1(t) & x \leq t \leq a, \\ f(t) & a \leq t \leq b, \\ g_2(t) & b \leq t \leq y, \\ t & y \leq t. \end{cases}$$

COROLLARY 3.5.1. *If $\mathcal{G}$ is an interpolating groupoid of germs on $\mathbf{R}$ that contains the identity, then $H(\mathcal{G})$ generates $\mathcal{G}$.*

LEMMA 3.6 (REDUCTION). *Let $\mathcal{G}$ be an interpolating groupoid of germs on $\mathbf{R}$ that contains the identity and is $s$-invariant with respect to the map $s(x) = x + r$. Let $\mathcal{G}'$ be the projection of $\mathcal{G}$ on $S_r$. Then there is a natural isomorphism from $N(\mathcal{G})/H(\mathcal{G})$ to $N(\mathcal{G}')/H(\mathcal{G}')$.*

PROOF: Let $h$ be in $N(\mathcal{G})$. We wish to find a function $g$ in the same coset of $H(\mathcal{G})$ in $N(\mathcal{G})$ so that $g$ is a lift of a homeomorphism of $S_r$. Note that if $g$ is increasing this means that we want $g$ to commute with $s$. However, if $g$ is decreasing, we want $g$ to satisfy $g = sgs$.

We know that $hsh^{-1}$ is in $H(\mathcal{G})$, as are $shs^{-1}h^{-1}$ and $shsh^{-1}$. Also, $shs^{-1}h^{-1}$ and $shsh^{-1}$ are increasing no matter what the behavior of $h$ is.

If $h$ is increasing, then $h(0) < h(r)$ and we let $a = h(0)$ and $c = h(r)$. Now $(shs^{-1}h^{-1})(h(r)) = h(0) + r$ and if we let $b = h(0)$ and $d = h(0) + r$, then we have $a < c$ and $b < d$. There is an $f$ in $H(\mathcal{G})$ whose germ at $h(0)$ is that of the identity and whose germ at $h(r)$ is that of $shs^{-1}h^{-1}$. Now $fh$ takes 0 to $h(0)$ with the germ of $h$ at 0 and takes $r$ to $h(0) + r$ with the germ of $shs^{-1}$ at $r$. Since the germ of $(fh)$ at $r$ is the germ of $s(fh)s^{-1}$ at $r$, we can build a function $g$ in $H(\mathcal{G})$ that agrees with $fh$ on $[0, r]$ and that commutes with $s$ on all of $\mathbf{R}$. Specifically, $g$ on $[ir, (i+1)r]$ is defined to be $s^i(fh)s^{-i}$. Slightly to the left of $ir$, the behavior of $g$ is that of $s^{i-1}(fh)s^{-i+1}$ with the relevant part of $fh$ being that in a small neighborhood of $r$. But here we have $fh = s(fh)s^{-1}$ so the behavior of $g$ to the left of $ir$ is also that of $s^i(fh)s^{-i}$ which agrees with the behavior of $g$ to the right of $ir$. Thus $g$ is $\mathcal{G}$-compatible at each $ir$. It is also $\mathcal{G}$-compatible in the interior of each $[ir, (i+1)r]$ since $\mathcal{G}$ is $s$-invariant.

If $h$ is decreasing, then $h(r) < h(0)$ and we let $a = h(r)$ and $c = h(0)$. Now $(shsh^{-1})(h(0)) = h(r) + r$ and if we let $b = h(r)$ and $d = h(r) + r$, then we have $a < c$ and $b < d$. There is an $f$ in $H(\mathcal{G})$ whose germ at $h(r)$ is that of the identity and whose germ at $h(0)$ is that of $shsh^{-1}$. Now $fh$ takes $r$ to $h(r)$ with the germ of $h$ at $r$ and takes 0 to $h(r) + r$ with the germ of $shs$ at 0. The rest of the argument is similar to the paragraph above except we define $g$ on $[-ir, -(i-1)r]$ to be $s^i f s^i$ and get a decreasing function that satisfies $sgs = g$.

In both cases, we get a function that defines a homeomorphism on $S_r = \mathbf{R}/r\mathbf{Z}$. This function is in the same coset as $h$ since it is a $\mathcal{G}$-rearrangement. We must show that another function $g'$ in the same coset as $h$ that has $sg's^{-1} = g'$ if $g'$ is increasing or $sg's = g'$ if $g'$ is decreasing has the projection of $g'$ to $S_r$ in the same coset as the projection of $g$. We know that $g^{-1}g'$ is in $H(\mathcal{G})$. Now $g'$ is increasing if and only if $g$ and $h$ are increasing. If $g$ and $g'$ are increasing, then $g^{-1}g' = sg^{-1}s^{-1}sg's^{-1} = sg^{-1}g's^{-1}$ and $g^{-1}g'$ projects to $S_r$ and shows that the projection of $g$ and $g'$ are in the same coset of $H(\mathcal{G}')$. If $h$ and $g'$ are decreasing, then $g^{-1}g' = s^{-1}g^{-1}s^{-1}sg's = s^{-1}g^{-1}gs$ and the same conclusion holds.

The map clearly defines a homomorphism. We want to show that the function is one to one and onto. It is onto by lifting representatives of $N(\mathcal{G}')/H(\mathcal{G}')$ to $\mathbf{R}$. Such lifts are in $N(\mathcal{G})$ since $\mathcal{G}$ is $s$-invariant. It is one to one, since if $h$ projects to a $\mathcal{G}'$-compatible homeomorphism of $S_r$, then $h$ is $\mathcal{G}$-compatible.

### 4. Normalizers and germs of Thompson's groups.

One of the goals in the paper is to show that non-PL homeomorphisms do not normalize the groups that we consider. We first reduce this to a study of normalizers of a small set of groups, and then show that we need only study conjugates of a small set of function. To do the first step, we need to compare the normalizers of the various groups. The next lemma starts this process.

LEMMA 4.1. *Let $G$ be a group or monoid for which $BPL_n(\mathbf{R}) \subseteq G \subseteq \widetilde{\overline{PL}_n}(\mathbf{R})$ or $BT_{n,r} \subseteq G \subseteq \widetilde{\overline{T}}_{n,r}$. Then $N(G) \subseteq N(BPL_n(\mathbf{R}))$ or $N(G) \subseteq N(BT_{n,r})$ whichever applies.*

PROOF: Conjugation preserves the properties of being invertible, being orientation preserving, having bounded support, and being the identity on an open interval.

We now know that the existence of a non-PL normalizer for any $G$ considered in Lemma 4.1 implies the existence of a non-PL element in either $N(BPL_n(\mathbf{R}))$ or $N(BT_{n,r})$. We wish to reduce the set of normalizers to consider to just the normalizers $N(BT_{n,n-1})$. The order of analysis is to first compare $N(BT_{n,r})$ to $N(BPL_n(\mathbf{R}))$, and then compare $N(BPL_n(\mathbf{R}))$ to $N(BT_{n,n-1})$.

We will work with groupoids of germs. Let $\mathcal{B}_{n,r}$ be the groupoid of germs generated by the group $BT_{n,r}$. Let $\mathcal{B}_n$ be the groupoid of germs generated by the group $BPL_n(\mathbf{R})$. From Lemma 2.5, we have $H(\mathcal{B}_n) = \ker(d_n)$. Let $\mathcal{P}_n$ be the groupoid of germs generated by the group $PL_n(\mathbf{R})$. We have $H(\mathcal{P}_n) = PL_n(\mathbf{R})$. The groupoid $\mathcal{P}_n$ was identified in Section 3 as $(\mathcal{A}_n(\mathbf{R}); \mathbf{Z}[\frac{1}{n}])$, but we will not need this fact and will use the simpler notation. Let $\widetilde{\mathcal{B}}_{n,r}$ be the lift of $\mathcal{B}_{n,r}$ to $\mathbf{R}$.

LEMMA 4.2. *We have the containments*

$$\mathcal{B}_n \subseteq \widetilde{\mathcal{B}}_{n,r} \subseteq \mathcal{P}_n$$

*and the equality*

$$\mathcal{B}_n = \widetilde{\mathcal{B}}_{n,n-1}.$$

PROOF: We get $\widetilde{\mathcal{B}}_{n,r} \subseteq \mathcal{P}_n$ because all lifts of elements of $BT_{n,r}$ are elements of $PL_n(\mathbf{R})$. The germs of the generators mentioned in Lemma 2.6 are contained in the germs of lifts of the generators mentioned in Lemma 2.7. The containment $\mathcal{B}_n \subseteq \widetilde{\mathcal{B}}_{n,r}$ now follows from Lemma 3.1.

The translation $x \mapsto x + (n-1)$ preserves $\Delta_n$ and so (using the notation of Lemma 2.7) $p_{n-1}^{-1}(p_{n-1}(\Delta_n)) = \Delta_n$. This and Lemma 2.7 imply that all lifts to $\mathbf{R}$ of elements of $BT_{n,n-1}$ preserve $\Delta_n$. Now Lemma 2.5 gives $\widetilde{\mathcal{B}}_{n,n-1} \subseteq \mathcal{B}_n$.

We now can compare $N(BT_{n,r})$ to $N(BPL_n(\mathbf{R}))$. Let $\widetilde{N}(BT_{n,r})$ be the set of all lifts to $\mathbf{R}$ of elements of $N(BT_{n,r})$.

LEMMA 4.3. *We have $\widetilde{N}(BT_{n,r}) \subseteq N(BPL_n(\mathbf{R}))$.*

PROOF: Let $\widetilde{BT}_{n,r}$ be the set of lifts to $\mathbf{R}$ of the elements of $BT_{n,r}$. We know that the elements of $\widetilde{N}(BT_{n,r})$ normalize $\widetilde{BT}_{n,r}$. However $\widetilde{BT}_{n,r}$ generates $\widetilde{\mathcal{B}}_{n,r}$, so $\widetilde{N}(BT_{n,r})$ normalizes $H(\widetilde{\mathcal{B}}_{n,r})$. We have

$$\mathcal{B}_n \subseteq \widetilde{\mathcal{B}}_{n,r} \subseteq \mathcal{P}_n$$

and thus

$$BPL_n(\mathbf{R}) \subseteq H(\mathcal{B}_n) \subseteq H(\widetilde{\mathcal{B}}_{n,r}) \subseteq H(\mathcal{P}_n) = PL_n(\mathbf{R}).$$

So elements of $\widetilde{N}(BT_{n,r})$ conjugate $BPL_n(\mathbf{R})$ into $PL_n(\mathbf{R})$. We now repeat the arguments in the proof of Lemma 4.1.

To compare the normalizers of $BPL_n(\mathbf{R})$ and $BT_{n,n-1}$ we need to know more about the groupoid of germs $\mathcal{B}_n$.

LEMMA 4.4.

  (i) $\mathcal{B}_n$ contains the identity.
  (ii) $\mathcal{B}_n$ is $s$-invariant where $s(x) = x + (n-1)$.
  (iii) $\mathcal{B}_n$ is interpolating.

PROOF: We get (i) and (ii) from the fact that the identity and $s$ are in $H(\mathcal{B}_n) = \ker(d_n)$. To see (iii) we let $a < c$ and $b < d$ where $(\mathcal{B}_n)_{a,b}$ and $(\mathcal{B}_n)_{c,d}$ are non-empty. Given $\alpha \in (\mathcal{B}_n)_{a,b}$ and $\beta \in (\mathcal{B}_n)_{c,d}$, we can get $a < a' < c' < c$ where $a'$ and $c'$ are in $\Delta_n$ and $[a, a']$ is in the domain of some representative $f$ of $\alpha$ and $[c', c]$ is in the domain of some representative $g$ of $\beta$. We know that $f$ and $g$ preserve $\Delta_n$. Now Lemma 2.3 allows us to fill the missing part on $[a', c']$.

COROLLARY 4.4.1. We have $N(\mathcal{B}_n)/H(\mathcal{B}_n) \simeq N(\mathcal{B}_{n,n-1})/H(\mathcal{B}_{n,n-1})$. In particular, if there is a non-PL normalizer of $BPL_n(\mathbf{R})$, then there is a non-PL normalizer of $BT_{n,n-1}$.

PROOF: The isomorphism follows from Lemmas 3.6, 4.2 and 4.4.

COROLLARY 4.4.2. Let $G$ be a group or monoid for which $BPL_n(\mathbf{R}) \subseteq G \subseteq \widetilde{PL}_n(\mathbf{R})$ or $BT_{n,r} \subseteq G \subseteq \overline{\widetilde{T}}_{n,r}$. If there is a non-PL normalizer of $G$, then there is a non-PL normalizer of $BT_{n,n-1}$.

PROOF: From Lemma 4.1, there is a non-PL element in either $N(BPL_n(\mathbf{R}))$ or $N(BT_{n,r})$. From Lemma 4.3, there is a non-PL element of $N(BPL_n(\mathbf{R}))$. Corollary 4.4.1 finishes the argument.

Let $\nu_n : S_{n-1} \to S_{n-1}$ be the degree $n$ map defined by $\nu_n(x) = nx$. It is an element of $\overline{T}_{n,n-1}$. It also preserves $p_{n-1}(\Delta_n)$ and so lies in $F(\mathcal{B}_{n,n-1})$.

LEMMA 4.5. *Let $G$ be a group or monoid for which $BPL_n(\mathbf{R}) \subseteq G \subseteq \widetilde{\overline{PL}}_n(\mathbf{R})$ or $BT_{n,r} \subseteq G \subseteq \widetilde{\overline{T}}_{n,r}$. If there is a non-PL normalizer of $G$, then there is a non-PL homeomorphism $h : S_{n-1} \to S_{n-1}$ for which $h(\mathbf{Z}[\frac{1}{n}]) = \mathbf{Z}[\frac{1}{n}]$ and for which $h\nu_n h^{-1}$ and $h^{-1}\nu_n h$ are in $\overline{T}_{n,n-1}$.*

PROOF: By Corollary 4.4.2, there is a non-PL normalizer $h$ of $BT_{n,n-1}$. Since $x \mapsto -x$ normalizes $BT_{n,n-1}$, we can assume that $h$ preserves orientation. Since the groups of germs $(\mathcal{B}_{n,n-1})_{x,x}$ at fixed points of functions in $BT_{n,n-1}$ are isomorphic to $\mathbf{Z} \times \mathbf{Z}$ when $x \in \mathbf{Z}[\frac{1}{n}]$, while at all other points these groups are isomorphic to $\mathbf{Z}$ or are trivial, we have that $h(\mathbf{Z}[\frac{1}{n}]) = \mathbf{Z}[\frac{1}{n}]$ from Lemma 3.3. We have $\nu_n \in F(\mathcal{B}_{n,n-1}) \subseteq \overline{T}_{n,n-1}$ and $BT_{n,n-1}$ generates $\mathcal{B}_{n,n-1}$. By Lemma 3.2, we get the last items of the statement.

## 5. Realizing automorphisms.

In this section we show that the automorphisms of the groups that we consider are all realized as conjugations by homeomorphisms. This follows from very general results in the theory of ordered groups. Our task is simply that of verifying the hypotheses. We also show that automorphisms of the monoids that we consider are realized as conjugations by homeomorphisms. The arguments simply reduce this to the group case.

We need hypotheses that define sufficient levels of transitivity for groups acting on $\mathbf{R}$ or $S^1$. For a group acting on $\mathbf{R}$, we use the notion of o-$k$-transitivity defined in Section 2. We need a different definition on $S^1$. Let $G$ be a group that acts on $S^1$ by homeomorphisms. We say that a closed interval (arc) is proper in $S^1$ if it and its complement have non-empty interior. We say that $G$ is approximately arc transitive if for all proper arcs $I$, $J_1$, $J_2$ in $S^1$, with $J_1 \subseteq \mathring{J}_2$, there is a $g \in G$ so that $J_1 \subseteq g(I) \subseteq J_2$. Note that if $G$ is approximately o-3-transitive (where "o" refers to preserving orientation), then $G$ is approximately arc transitive. Recall that a homeomorphism of $S^1$ has bounded support if it fixes some open interval. Most of the work in this section is taken care of in the following theorem from [18].

THEOREM 5.1 (MCCLEARY-RUBIN). *Let $G$ act on $\mathbf{R}$ or $S^1$ by homeomorphisms. Assume that $G$ contains a non-identity element of bounded support. If $G$ acts on $\mathbf{R}$, then assume that the action is approximately o-2-transitive, and if $G$ acts on $S^1$, then assume that the action is approximately arc transitive. Then for each automorphism $\alpha$ of $G$, there is a unique self homeomorphism $h$ of $\mathbf{R}$ or $S^1$ (whichever is appropriate) so that $\alpha(f) = hfh^{-1}$ for every $f \in G$.*

We will use this to demonstrate provision (i) of Theorem 1. Specifically, we will show:

THEOREM 1. *Let $G$ be a group or monoid for which $BPL_n(\mathbf{R}) \subseteq G \subseteq \widetilde{PL_n}(\mathbf{R})$ or $BT_{n,r} \subseteq G \subseteq \widetilde{T}_{n,r}$. If $G$ is not a group, then assume $G$ is determined by its invertible elements. Then*

  (i) *$\Phi : N(G) \to Aut(G)$ is an isomorphism.*

The first step is to show:

THEOREM 5.2. *Let $G$ be a group for which $BPL_n(\mathbf{R}) \subseteq G \subseteq \widetilde{PL_n}(\mathbf{R})$ or $BT_{n,r} \subseteq G \subseteq \widetilde{T}_{n,r}$. Then $\Phi : N(G) \to Aut(G)$ is an isomorphism.*

PROOF: That $G$ satisfies the hypotheses of Theorem 5.1 follows from the results in Section 2. The existence and uniqueness statements in Theorem 5.1 show that $\Phi$ is one to one and onto.

We use Theorem 5.2 in the proof of the next theorem which handles the case in which $G$ is a monoid.

THEOREM 5.3. *Let $M$ be a monoid for which $BPL_n(\mathbf{R}) \subseteq M \subseteq \widetilde{PL_n}(\mathbf{R})$ or $BT_{n,r} \subseteq M \subseteq \widetilde{T}_{n,r}$ and assume $M$ is determined by its invertible elements. Then $\Phi : N(M) \to Aut(M)$ is an isomorphism.*

PROOF: Note that the group of invertible elements in $M$ satisfies the hypotheses of Theorem 5.2. We will refer to this group as $G(M)$. If $\phi$ is an automorphism of $M$, then $\phi|_{G(M)}$ is an automorphism of $G(M)$. By Theorem 5.2, there is an $h \in \text{Homeo}(X)$ so that $\phi|_{G(M)}$ is realized as conjugation by $h$. By Corollary 3.2.1, $\Phi(h)$ is an automorphism of $M$. If we form the composition $(\Phi(h))^{-1} \circ \phi$, then we have an automorphism of $M$ that is the identity on $G(M)$. We are done if we show that $(\Phi(h))^{-1} \circ \phi$ is the identity on all of $M$. We thus reduce the problem and assume that $\phi$ is an automorphism on $M$ that is the identity on $G(M)$. Our task is to show that $\phi$ is the identity on $M$.

We focus on $P$, the set of restrictions of the identity to open subsets of $X$. These are all in $M$ because $M$ is determined by its invertible elements, among which is the identity of $M$. All the elements of $P$ are idempotents of $M$, $\alpha \in P$ implies $\alpha\alpha = \alpha$, but they are not all the idempotents. We wish to distinguish $P$ from the other idempotents. For $\alpha \in P$, if $f\alpha = \alpha$, then $f$ must the identity on the domain of $\alpha$ and is arbitrary off the domain. If $f$ is invertible, then it must preserve the complement of the domain of $\alpha$. This will give $\alpha f = \alpha$. Thus every invertible left identity of $\alpha$ is a right identity of $\alpha$. Now suppose that $\beta$ is an idempotent that has $\beta(x) \neq x$ for some $x$. There is an open set $U$ about $x$ in the domain of $\beta$ for which $\beta(U) \cap U = \emptyset$. By further restricting $U$, we can assume that $\beta$ is one to one on $U$. If some $y \in U$ is in the range of $\beta$, then $y = \beta(z)$ and

$$y \neq \beta(y) = \beta(\beta(z)) = \beta(z) = y$$

so $U$ is disjoint from the range of $\beta$. Using the transitivity properties of $G(M)$, we can build an invertible $g$ in $M$ that is the identity on the range of $\beta$ and takes $x$ to some point $w \in U$ with $x \neq w$. Since $g$ is the identity on the range of $\beta$ it is a left identity for $\beta$. But if $g$ is also a right identity for $\beta$, then we get $\beta(x) = \beta(g(x)) = \beta(w)$. With both $x$ and $w$ in $U$ on which $\beta$ is one to one, we get $x = w$. Thus $g$ cannot be a right identity for $\beta$. Since we have characterized the elements of $P$ using the operation in $M$, we have that $\phi(P) = (P)$.

We want to argue that $\phi$ is the identity on $P$. An element of $P$ is identified by its domain. For $\alpha \in P$, let $D(\alpha)$ be the closure of the domain of $\alpha$. Assume that $D(\alpha)$ and $D(\phi(\alpha))$ are different. Assume first that $D(\phi(\alpha))$ has points that are not in $D(\alpha)$. Thus the domain of $\phi(\alpha)$ has a point $x$ not in $D(\alpha)$ and there is an open set $U$ about $x$ in the domain of $\phi(\alpha)$ that is disjoint from $D(\alpha)$. Thus we can build an invertible element $g$ in $M$ that is the identity on $D(\alpha)$ but has $g(x) \neq x$. Recall that $\alpha$ and $\phi(\alpha)$ are the identities on their domains. Since $g$ is invertible, it is fixed by $\phi$ and since $g$ is the identity on $D(\alpha)$ we get $g\alpha = \alpha$. Thus $\phi(\alpha) = \phi(g\alpha) = \phi(g)\phi(\alpha) = g\phi(\alpha)$. This says that $x = \phi(\alpha)(x) = g\phi(\alpha)(x) = g(x) \neq x$. Thus $D(\phi(\alpha)) \subseteq D(\alpha)$. By working with $\phi^{-1}$ we also get $D(\alpha) \subseteq D(\phi(\alpha))$ and so we must have $D(\phi(\alpha)) = D(\alpha)$.

Now let $d(\alpha)$ be the domain of $\alpha \in P$. Assume that $d(\alpha)$ contains a point $x$ that is not in $d(\phi(\alpha))$. By the transitivity properties of $G(M)$ and since $d(\alpha)$ is an open set, we can build an invertible element $g$ of $M$ whose support lies in $d(\alpha)$ and for which $g(x) \neq x$. Now $\phi(\alpha)$ is not defined on $x$, so $\phi(\alpha)g^{-1}$ cannot be defined on $g(x)$. Because the support of $g$ lies in $d(\alpha)$, so does the support of $g^{-1}$, and we have $g\alpha g^{-1} = \alpha$. Since $g$ and $g^{-1}$ are invertible, we have $\phi(g) = g$, $\phi(g^{-1}) = g^{-1}$ and

$$\phi(\alpha) = \phi(g\alpha g^{-1}) = g\phi(\alpha)g^{-1}.$$

But the right hand expression is not defined on $g(x)$ so neither is $\phi(\alpha)$. Since we can build sufficiently varied invertible elements $g$ in $M$ to carry $x$ to a dense subset of $d(\alpha)$, we get that $\phi(\alpha)$ is not defined on a dense subset of $d(\alpha)$. The domain of $\phi(\alpha)$ is an open set, so we now know that $d(\phi(\alpha))$ is disjoint from $d(\alpha)$. But we showed that $d(\alpha)$ and $d(\phi(\alpha))$ have identical closures, so we must have $d(\alpha) \subseteq d(\phi(\alpha))$. Again, by working with $\phi^{-1}$ we get that $d(\phi(\alpha)) \subseteq d(\alpha)$. We have shown that $\phi$ preserves the domains of the elements of $P$, and from this it follows that $\phi$ fixes the elements of $P$.

We can use elements of $P$ to identify domains. We have an order on elements of $P$ given by subset on the domains. This is captured by the monoid operation by saying that $f$ and $g$ in $P$ are related by $f \leq g$ if and only if $gf = fg = f$. We can now identify the domain of an element $h$ in $M$ as the minimum $f$ in $P$ that satisfies $hf = h$. Now that we know that $\phi$ fixes elements of $P$, we know that $\phi$ preserves domains in that for $h$ in $M$, the domain of $\phi(h)$ equals the domain of $h$.

We can use subsets of $P$ to identify values of functions. Let $x$ be in $X$ and let $U_x = \{f \in P \mid f(x) \text{ is defined}\}$. Let $h$ be in $M$ with $x$ in the domain of $h$. The statement that $h(x) = y$ is equivalent to the statement that for all $f \in U_x$ and all $g \in U_y$, $ghf$ is not the empty function. Since the empty function is a zero (it is absorbing under composition on either side with all functions) it is determined by the operations. Since $\phi$ fixes the empty function and all the elements of $P$, it preserves the statement $h(x) = y$ in that $\phi(h)(x) = y$ for all $x \in X$ and all $h \in M$. Thus $\phi$ is the identity on $M$ which is what we wanted to show.

## PART III. PL NORMALIZERS

### 6. Proof of Theorems 1 and 2 from Theorem 3.

It is the goal of this section to prove Theorems 1 and 2 and Corollary 2.1 assuming Theorem 3. These implications are true even when $n$ is not restricted to 2. It is only in the proof in Part IV of Theorem 3 that we need to assume $n = 2$. We restate Theorems 1, 2, Corollary 2.1 and Theorem 3 for arbitrary $n$. We call these "Statements" since only the implications between them are known to be true. We also show the "equivalence" of Theorem 2 and Corollary 2.1 by showing that Statement 2.1′ implies Statement 2′ when $n = 2$.

STATEMENT 1′. *Let $G$ be a group or monoid for which $BPL_n(\mathbf{R}) \subseteq G \subseteq \widetilde{\overline{PL}}_n(\mathbf{R})$ or $BT_{n,r} \subseteq G \subseteq \widetilde{\overline{T}}_{n,r}$. If $G$ is not a group, then assume $G$ is determined by its invertible elements. Then*

  (i) *$\Phi : N(G) \to Aut(G)$ is an isomorphism,*
  (ii) *$N(G) \subseteq \widetilde{PL}_n(\mathbf{R})$ or $N(G) \subseteq \widetilde{T}_{n,r}$ whichever applies,*
  (iii) *the containment is equality if $G$ is one of $BPL_n(\mathbf{R})$, $PL_n(\mathbf{R})$, $\widetilde{PL}_n(\mathbf{R})$, $\overline{PL}_n(\mathbf{R})$, $\widetilde{\overline{PL}}_n(\mathbf{R})$, $BT_{n,r}$, $T_{n,r}$, $\widetilde{T}_{n,r}$, $\overline{T}_{n,r}$, or $\widetilde{\overline{T}}_{n,r}$, and*
  (iv) *if $N_o(F_n)$ represents the index 2 subgroup of $N(F_n)$ of orientation preserving elements, then there is an exact sequence*

$$1 \to F_n \to N_o(F_n) \to T_{n,n-1} \times T_{n,n-1} \to \mathbf{Z}_{n-1} \to 1.$$

STATEMENT 2′. *Let $h : S_{n-1} \to S_{n-1}$ be an orientation preserving homeomorphism for which $h(\mathbf{Z}[\frac{1}{n}]) = \mathbf{Z}[\frac{1}{n}]$, and assume that $h\nu_n h^{-1}$ and $h^{-1}\nu_n h$ are both in $\overline{T}_{n,n-1}$. Then $h \in T_{n,n-1}$.*

STATEMENT 2.1′. *Let $h : S_{n-1} \to S_{n-1}$ be an orientation preserving homeomorphism that fixes the images of $\mathbf{Z}$ in $S_{n-1}$, and assume that $h\nu_n h^{-1}$ and $h^{-1}\nu_n h$ are both in $\overline{T}_{n,n-1}$. Then $h \in T_{n,n-1}$.*

STATEMENT 3′. *Let $h : S_{n-1} \to S_{n-1}$ be an orientation preserving homeomorphism for which $h(\mathbf{Z}[\frac{1}{n}]) = \mathbf{Z}[\frac{1}{n}]$, and assume that $h\nu_n h^{-1}$ and $h^{-1}\nu_n h$ are both in $\overline{T}_{n,n-1}$. Then $h$ is PL.*

For the rest of this section, we assume the truth of Statement 3′. From Lemma 4.5, we know the following.

STATEMENT 6.1. *Let $G$ be a group or monoid for which $BPL_n(\mathbf{R}) \subseteq G \subseteq \widetilde{PL}_n(\mathbf{R})$ or $BT_{n,r} \subseteq G \subseteq \widetilde{\overline{T}}_{n,r}$. Then every element in the normalizer of $G$ is PL.*

PROOF OF STATEMENT 1′ FROM STATEMENT 3′: Part (i) is proven in Section 5. From Lemmas 4.1 and 4.3, we know that if we prove that $N(BPL_n(\mathbf{R})) \subseteq \widetilde{PL}_n(\mathbf{R})$, then we can conclude that $N(G) \subseteq \widetilde{PL}_n(\mathbf{R})$ or $N(G) \subseteq \widetilde{T}_{n,r}$ for the objects $G$ covered in the lemmas. This will prove Part (ii) of Statement 1′. Thus we must prove:

STATEMENT 6.2. $N(BPL_n(\mathbf{R})) \subseteq \widetilde{PL}_n(\mathbf{R})$.

PROOF ASSUMING STATEMENT 3′: Let $h$ be an element of $N(BPL_n(\mathbf{R}))$. As argued in the proof of Lemma 4.5, we know that $h(\mathbf{Z}[\frac{1}{n}]) = \mathbf{Z}[\frac{1}{n}]$. Given any interval on which $h$ is affine, there are a pair of elements of $\mathbf{Z}[\frac{1}{n}]$ in the interval whose distance apart is in integral power of $n$. Since the images of these two points are in $\mathbf{Z}[\frac{1}{n}]$, the distance between the images is in $\mathbf{Z}[\frac{1}{n}]$, and the slope of $h$ on the interval must have an integral power of $n$ for the denominator. The discussion also applies to $h^{-1}$, and it follows that the slope is an integral power of $n$. This completes the proof of the this statement and Part (ii) of Statement 1′.

The proof of Part (iii) of Statement 1′ is one of checking that $\widetilde{PL}_n(\mathbf{R})$ and $\widetilde{T}_{n,r}$ normalize the groups listed. This is quite easy and is left to the reader.

We consider (iv) of Statement 1′. We know that $N(F_n) \subseteq \widetilde{PL}_n(\mathbf{R})$. Let $f$ be in $N_o(F_n)$, the orientation preserving elements of $N(F_n)$. The germ of $F_n$ at $-\infty$ is isomorphic to $\mathbf{Z}$ generated by translation to the right or left by $(n-1)$. A normalizer must preserve the germ at $-\infty$ and so must take a generator to a generator. (This is identical to the concept discussed in Lemma 3.3.) An orientation preserving conjugator will take a right translation to a right translation. Thus near $-\infty$, an orientation preserving normalizer commutes with translation to the right by $(n-1)$. This forces $f(x+n-1) = f(x)+(n-1)$ near $-\infty$. Similarly, we get $f(x + n - 1) = f(x) + (n - 1)$ near $+\infty$. Conversely, any $f$ in $PL_n(\mathbf{R})$ that satsifies $f(x+n-1) = f(x)+(n-1)$ near $\pm\infty$ conjugates elements of $F_n$ into $F_n$. Thus we have characterized the elements of $N_o(F_n)$. Now by restricting an $f \in N_o(F_n)$ sufficiently far to the left, we get a well defined element $f_-$ of $T_{n,n-1}$. Similarly, the restriction of $f$ near $+\infty$ gives a well defined element $f_+$ of $T_{n,n-1}$. The function $f \mapsto (f_-, f_+)$ is a homomorphism into $T_{n,n-1} \times T_{n,n-1}$. The kernel consists of all elements of $F_n$.

We use the homomorphism $d_n$ of Section 2 to understand the image of the homomorphism above. Each element $g$ of $T_{n,n-1}$ lifts to various elements of $PL_n(\mathbf{R})$ that all have the same value of $d_n$. This gives a well defined meaning to $d_n(g)$. Easy applications of Lemmas 2.2 and 2.3 show that two elements $g_1$ and $g_2$ of $T_{n,n-1}$ will have the lift of $g_1$ near $-\infty$ and the lift of $g_2$ near $+\infty$ connectable into one function in $PL_n(\mathbf{R})$ if and only if $d_n(g_1) = d_n(g_2)$. Thus the image of the homomorphism is the kernel of $(g_1, g_2) \mapsto d_n(g_1) - d_n(g_2)$. This completes the proof of Statement $1'$.

PROOF OF STATEMENT $2'$ FROM STATEMENT $3'$: From Statement $1'$, Part (ii), we are done when we show the truth of the following lemma.

LEMMA 6.3. *Let $h : S_{n-1} \to S_{n-1}$ be an orientation preserving homeomorphism for which $h(\mathbf{Z}[\frac{1}{n}]) = \mathbf{Z}[\frac{1}{n}]$, and assume that $h\nu_n h^{-1}$ and $h^{-1}\nu_n h$ are both in $\overline{T}_{n,n-1}$. Then $h$ is in $N(BT_{n,n-1})$.*

REMARK: This lemma does not depend on Statement $3'$, and is true for all $n$ as stated. This lemma will also be important in the proof of Theorem 3 in Part IV.

PROOF: Let $\mathcal{V}_n$ be the groupoid of germs generated by the function $\nu_n$. From Lemma 3.2, we know that $h$ normalizes $H(\mathcal{V}_n; \mathbf{Z}[\frac{1}{n}])$. Recall that $\mathcal{B}_{n,n-1}$ is the groupoid of germs generated by $BT_{n,n-1}$. We have $\mathcal{B}_{n,n-1} \subseteq \mathcal{T}_{n,n-1}$ where $\mathcal{T}_{n,n-1}$ is the groupoid of germs generated by $T_{n,n-1}$. Thus $H(\mathcal{B}_{n,n-1}) \subseteq H(\mathcal{T}_{n,n-1}) = T_{n,n-1}$. If we show that $(\mathcal{V}_n; \mathbf{Z}[\frac{1}{n}]) = \mathcal{B}_{n,n-1}$, then we have that $h$ normalizes a group $G$ with $BT_{n,n-1} \subseteq G \subseteq T_{n,n-1}$, and by Lemma 4.1, we are done. Thus we must show the following.

LEMMA 6.4. $(\mathcal{V}_n; \mathbf{Z}[\frac{1}{n}]) = \mathcal{B}_{n,n-1}$.

PROOF: Since translation by $n-1$ preserves $\Delta_n$, we have that the ring homomorphism $\phi_n : \mathbf{Z}[\frac{1}{n}] \to \mathbf{Z}_{n-1}$ makes sense on $S_{n-1}$. Since $\nu_n$ is just a shift of the expression $x = a_j a_{j-1} \ldots a_0.a_{-1}a_{-2}\ldots$ as in (3-2), it commutes with $\phi_n$ and preserves $p_{n-1}(\Delta_n)$. Combining this with Lemma 2.7 gives $\mathcal{V}_n \subseteq \mathcal{B}_{n,n-1}$. It follows from Lemma 2.5 that arbitrary changes of slope of elements of $BT_{n,n-1}$ can occur at elements of $\mathbf{Z}[\frac{1}{n}]$. This puts $(\mathcal{V}_n; \mathbf{Z}[\frac{1}{n}]) \subseteq \mathcal{B}_{n,n-1}$.

Since $(\mathcal{V}_n; \mathbf{Z}[\frac{1}{n}])$ allows arbitrary changes of slope at elements of $\mathbf{Z}[\frac{1}{n}]$, we need only show that all linear parts of elements of $BT_{n,n-1}$ are in $\mathcal{V}_n$. If we are presented with a linear part of an element $f$ of $BT_{n,n-1}$, then it is the restriction to an interval $I$ of an affine map $x \mapsto xn^k + d$ where $d \in \Delta_n$. We can assume that the interval $I$ has an element $a$ of $\mathbf{Z}[\frac{1}{n}]$ as a left endpoint. If $f(a) = b$, then we must show that $(\nu_n)^j(a) = (\nu_n)^k(b)$ for some integers $j$ and $k$ so that $(\nu_n)^k \circ f = (\nu_n)^j$ on $I$. This will make $f|_I$ the composition of two germs in $\mathcal{V}_n$. However, each of $a$ and $b$ have finite expression of the type in (3-2). The action of $\nu_n$ is to shift the expressions to the left and reduce modulo $(n-1)$. Under powers of $\nu_n$, the expressions for $a$ and $b$ can be reduced to single integers which represent $\phi_n(a)$ and $\phi_n(b)$. Since $f$ is in $BT_{n,n-1}$ which preserves $p_{n-1}(\Delta_n)$, and lifts of $f$ commute with $\phi_n$, we have that $\phi_n(a) = \phi_n(b)$ and thus have a common point in the forward orbit of $\nu_n$. This completes the proof.

REMARK: The ideas behind Lemma 6.3, Corollary 4.4.2 and Lemma 4.5 motivated the approach taken in this paper. Lemma 6.3, Corollary 4.4.2 and Lemma 4.5 make it clear that if Statement 1′ were to fail, then it would be possible to find examples of this failure by investigating conjugates of the maps $\nu_n$. Since the maps $\nu_n$ are expanding, the conjugators are determined by the conjugates.

We can give a separately worded result concerning rigidity. It is simply a translation of parts of Theorem 1 into other terminology. There is a corresponding statement for arbitrary $n$ that would be true if Statement 3′ is true. Recall that $\mathcal{B}_{n,r}$, $\mathcal{T}_{n,r}$, $\mathcal{B}_n$ and $\mathcal{P}_n$ were defined as the groupoids of germs generated by $BT_{n,r}$, $T_{n,r}$, $BPL_n(\mathbf{R})$ and $PL_n(\mathbf{R})$ respectively. Recall from Section 3 that if $\mathcal{G}$ is a groupoid of germs on $X$, then a homeomorphism $h: X \to X$ preserves $\mathcal{G}$ under conjugation if $hfh^{-1}$ is in $F(\mathcal{G})$ for all $f \in F(\mathcal{G})$. We say that $h$ normalizes $\mathcal{G}$ if $h$ and $h^{-1}$ preserve $\mathcal{G}$ under conjugation.

THEOREM 6.5. *If an orientation preserving homeomorphism $h : \mathbf{R} \to \mathbf{R}$ normalizes $\mathcal{B}_2$ or $\mathcal{P}_2$, then $h$ is $\mathcal{P}_2$-compatible. If an orientation preserving homeomorphism $h : S_1 \to S_1$ normalizes $\mathcal{B}_{2,1}$ or $\mathcal{T}_{2,1}$, then $h$ is $\mathcal{T}_{2,1}$-compatible.*

PROOF OF STATEMENT 2.1′ FROM STATEMENT 2′: We must show $h(\mathbf{Z}[\frac{1}{n}]) = \mathbf{Z}[\frac{1}{n}]$. The fixed points of $\nu_n$ are the integer points in $S_{n-1}$. The preimages of the fixed points under powers of $\nu_n$ are exactly the points of $\mathbf{Z}[\frac{1}{n}]$ in $S_{n-1}$. By hypothesis, $h$ fixes the fixed points of $\nu_n$. Thus the conjugates of $\nu_n$ by $h$ and $h^{-1}$ also have their fixed points at the integer points of $S_{n-1}$ and $h$ and $h^{-1}$ carry $\mathbf{Z}[\frac{1}{n}]$ (the preimages of the fixed points of powers $\nu_n$) into the preimages of the fixed points of the powers of the conjugates. Since the conjugates of $\nu_n$ by $h$ and $h^{-1}$ are in $T_{n,n-1}$, we know that they preserve $\mathbf{Z}[\frac{1}{n}]$ and $S_{n-1} - \mathbf{Z}[\frac{1}{n}]$ so that the preimages of the fixed points of powers of the conjugates are all in $\mathbf{Z}[\frac{1}{n}]$. Thus both $h$ and $h^{-1}$ carry $\mathbf{Z}[\frac{1}{n}]$ to $\mathbf{Z}[\frac{1}{n}]$.

PROOF OF STATEMENT 2′ FROM STATEMENT 2.1′ WHEN $n = 2$: Let $h$ be as in the hypothesis of Statement 2′. We are conjugating $\nu_2$ by a homeomorphism $h$ and its inverse. The unique fixed point of $\nu_2$ is the point 0 in $S_1$. Since $h(\mathbf{Z}[\frac{1}{2}]) = \mathbf{Z}[\frac{1}{2}]$, we know that the point 0 is taken to an element $x$ of $\mathbf{Z}[\frac{1}{2}]$. We know that rotation $\rho$ by $x$ is in $T_{2,1}$ and thus normalizes $\overline{T}_{2,1}$. We are done if we show that $\rho^{-1}h$ is in $T_{2,1}$. But this is true since $\rho^{-1}h$ satisfies the hypotheses of Statement 2.1′.

## PART IV. NON-PL NORMALIZERS — PROOF OF THEOREM 3

In this part we will prove:

THEOREM 3. *Let $h : S_1 \to S_1$ be an orientation preserving homeomorphism for which $h(\mathbf{Z}[\frac{1}{2}]) = \mathbf{Z}[\frac{1}{2}]$, and assume that $h\nu_2 h^{-1}$ and $h^{-1}\nu_2 h$ are both in $\overline{T}_{2,1}$. Then $h$ is PL.*

We had hoped to prove:

STATEMENT 3′. *Let $h : S_{n-1} \to S_{n-1}$ be an orientation preserving homeomorphism for which $h(\mathbf{Z}[\frac{1}{n}]) = \mathbf{Z}[\frac{1}{n}]$, and assume that $h\nu_n h^{-1}$ and $h^{-1}\nu_n h$ are both in $\overline{T}_{n,n-1}$. Then $h$ is PL.*

All of the sections in this part except Section 11 apply to Statement 3′. Thus except for Section 11, we will give arguments that apply to all $n$ and not just $n = 2$.

Our outline will be to extract as much information as possible about $h$ from the fact that $h\nu_n h^{-1}$ is in $\overline{T}_{n,n-1}$. Then, when $n = 2$, we will prove that $h$ must be PL if $h^{-1}\nu_n h$ is PL. It might seem that we will only use that $h^{-1}\nu_n h$ is PL. However, we will also use Lemma 6.3 which uses the full hypothesis. We do not know if we can weaken the hypotheses of Theorem 3 so that we only assume that $h^{-1}\nu_n h$ is PL. See Example 2 and remarks in Section 13. The use of Lemma 6.3 will be crucial in our argument, and we point out that its conclusion allows us to assume throughout this part that $h$ is a normalizer of $BT_{n,n-1}$. Since any rotation by an element of $\mathbf{Z}[\frac{1}{n}]$ is a normalizer of $BT_{n,n-1}$, we can compose $h$ by a rotation through $-h(0)$ to get an element satisfying the hypotheses that fixes 0 and that is PL if and only if $h$ is. Thus throughout the rest of this part, we assume $h(0) = 0$.

For the remainder of this part we will let $g = h\nu_n h^{-1}$.

## 7. Markov partitions.

We can build a very nice Markov partition for $g$. It is derived from a Markov partition for $\nu_n$. A Markov partition for the map $\nu_n$ with $(n-1)n^k$ intervals can be obtained by dividing $S_{n-1}$ into the intervals $[in^{-k}, (i+1)n^{-k}]$ for all integers $i$ with $0 \leq i < (n-1)n^k$. The action of the map $\nu_n$ is to take each interval affinely onto the union of $n$ consecutive intervals. We call this the standard Markov partition for $\nu_n$ of level $k$. The endpoints of these intervals are all in $\mathbf{Z}[\frac{1}{n}]$. Given finitely many points in $\mathbf{Z}[\frac{1}{n}]$ we can realize them as a subset of the endpoints of such a partition by taking $k$ large enough.

We can now apply $h$ to the endpoints of these intervals and get a sequence of points in $S_{n-1}$ that divide $S_{n-1}$ into $(n-1)n^k$ intervals that form a Markov partition for $g$. We invent notation for these intervals by defining

(7-1)
$$I_i^k = h([in^{-k}, (i+1)n^{-k}]).$$

We treat the subscripts cyclically modulo $(n-1)n^k$. We let $P^k$ denote the Markov partition itself. That is, $P^k$ is the collection of intervals $\{I_i^k \mid 0 \leq i < (n-1)n^k\}$. Properties of the standard Markov partition of level $k$ for $\nu_n$ are carried over by $h$ to the partition $P^k$.

The action of $\nu_n$ on $[in^{-k}, (i+1)n^{-k}]$ is to take it affinely onto the union of $[jn^{-k}, (j+1)n^{-k}]$ for $ni \leq j < n(i+1)$. The action of $g$ on the intervals of $P^k$ is to take $I_i^k$ homeomorphically onto the union of $I_j^k$ for $ni \leq j < n(i+1)$. We argue that this action fails to be affine for only finitely many of these intervals. We know that $h(\mathbf{Z}[\frac{1}{n}]) = \mathbf{Z}[\frac{1}{n}]$. With $g \in \overline{T}_{n,n-1}$, we have the set $D$ of break points of $g$ (the set $x$ on which $g'(x)$ fails to exist) is a finite subset of $\mathbf{Z}[\frac{1}{n}]$. The set $h^{-1}(D)$ is a finite set in $\mathbf{Z}[\frac{1}{n}]$ and is in the set of endpoints for the standard Markov partition for $\nu_n$ of level $k$ for some finite $k$. We let $K$ be the smallest value of $k$ for which this applies and we call this the stable level for $g$. We will need to refer to $K$ repeatedly in what follows and the reason for the word stable will become clear later. We now have that the set $D$ of breakpoints of $g$ is contained in the set of endpoints of the intervals of $P^K$. This makes the action of $g$ affine on each interval in $P^K$. The endpoints of $P^K$ are among the endpoints of $P^k$ for all $k \geq K$. Thus the action of $g$ affine on each interval in $P^k$ for all $k \geq K$.

We have Markov partitions $P^k$ for $g$ for every positive integer $k$. Each is obtained from the corresponding stardard Markov partition of level $k$ for $\nu_n$ by applying $h$. We can also describe how each $P^{k+1}$ is obtained from $P^k$. Each interval $I_i^k$ is subdivided into $n$ intervals $I_j^{k+1}$ for $ni \leq j < n(i+1)$ by adding as new endpoints those point in the interior of $I_i^k$ that are carried to the endpoints of $I_j^k$ for $ni \leq j < n(i+1)$. That is, we are pulling back the structure $P^k$ locally by the action of $g^{-1}$. Since $g$ is a local homeomorphism, this makes sense. We can thus refer to $P^{k+1}$ as the partition derived from $P^k$.

When $k \geq K$, the description gets simpler. Now, $I_i^k$ is subdived into intervals $I_j^{k+1}$, $ni \leq j < n(i+1)$, whose succession of lengths in increasing $j$ are in proportion to the succession of lenghts in increasing $j$ of the intervals $I_j^k$, $ni \leq j < n(i+1)$. This is because with $k > K$, the action of $g$ on each of the intervals of $P^k$ is affine.

Not only does $g$ carry intervals of $P^k$ onto unions of intervals of $P^k$, but we can describe how $g$ carries intervals of one partition to intervals of another. We note that $g$ carries $I_j^{k+1}$ of $P^{k+1}$ homeomorphically onto $I_j^k$ of $P^k$. This action is affine if $k+1 \geq K$. Inductively, we get

(7-2) $\qquad g^m$ carries $I_j^{k+m}$ of $P^{k+m}$ homeomorphically onto $I_j^k$ of $P^k$.

Again, the action is affine if $k+1 \geq K$ (not just $k+m \geq K$). Note that the subscripts of the intervals of $P^k$ are treated cyclically modulo $(n-1)n^k$ and the subscripts of the intervals $P^{k+m}$ are treated cyclically modulo $(n-1)n^{k+m}$.

We develop more notation that will be needed later. Each $I_i^k$ in $P^k$ has two endpoints, and we label them so that
$$I_i^k = [x_i^k, x_{i+1}^k].$$

We use $VP^k$ to denote the endpoints (vertices) of the partition $P^k$. Since endpoints of the partition $P^k$ are also endpoints in all partitions $P^{k+m}$, we have $VP^k \subseteq VP^{k+m}$ and we end up with multiple labels for the same points. Thus $x_i^k = x_{in^m}^{k+m}$. Note that since $I_i^k = h([in^{-k}, (i+1)n^{-k}])$, we have

(7-3) $\qquad\qquad\qquad x_i^k = h(in^{-k}) = h(i/n^k).$

We know that every $x \in \mathbf{Z}[\frac{1}{n}]$ is of such form and if $x \in \mathbf{Z}[\frac{1}{n}]$ has $x = h(i/n^k)$ so that $n \nmid i$, then $k$ is the least integer for which $x \in VP^k$. We call this the natural level of $x$ and denote it $\lambda(x)$. We can derive some information about $\lambda(x)$.

Let $x = h(i/n^k)$ so that $n \nmid i$ and $\lambda(x) = k$. Since $i/n^k$ is in $S_{n-1} = \mathbf{R}/(n-1)\mathbf{Z}$, we can assume that $0 \leq i/n^k < (n-1)$. If $i/n^k$ is expanded in base $n$ as given in (3-2), then there will be no more than one non-zero numeral to the left of the $n$-ary point, and the number of non-zero numerals to the right of the $n$-ary point is exactly $k = \lambda(x)$. If $k = 0$, then $i/n^k = i$ is one of the integers $0 \leq i < (n-1)$. The action of $\nu_n$ on the expansion (3-2) is to shift the expression one position to the left and then reduce modulo $(n-1)$. If $k = 0$, then $i/n^k = i$ is a fixed point of $\nu_n$ and $\lambda(g(x)) = \lambda(x)$. Otherwise $\nu_n(i/n^k)$ has one less numeral to the right of the $n$-ary point than $i/n^k$ and $\lambda(g(x)) = \lambda(x) - 1$. In general

(7-4) $\qquad\qquad\qquad \lambda(g^m(x))$ is the larger of $\lambda(x) - m$ or $0$.

We can write down the orbits of an element of $VP^k$. If $x \in \mathbf{Z}[\frac{1}{n}]$ has $x = h(i/n^k) = x_i^k$, then the action of $g$ on $x$ is to take it to $g(x) = h(ni/n^k) = h(i/n^{k-1}) = x_{ni}^k = x_i^{k-1}$. The action of higher powers of $g$ has

(7-5) $\qquad\qquad\qquad g^m(x_i^{k+m}) = x_{in^m}^{k+m} = x_i^k.$

Again the subscripts are treated cyclically modulo a base that varies with the superscript. The natural level of an $x$ changes by 1 with each appliation of $g$ and eventually, the image of an $x$ under powers of $g$ reach level zero. The points of level zero are the images under $h$ of the integers in $S_{n-1}$. The integers in $S_{n-1}$ are the fixed points of $\nu_n$, so their images under $h$ are the fixed points of $g$. Thus each $x = x_i^k$ has a fixed point of $g$ in its forward orbit under $g$. The value of this fixed point is $x_{i'}^0 = h(i')$ where $i' \equiv i \mod (n-1)$. Thus the forward orbits under $g$ of elements of $\mathbf{Z}[\frac{1}{n}]$ arrange themselves into $(n-1)$ different classes depending at which fixed point the forward orbits ends. We call these classes the extended orbits of $g$. Recalling the ring homomorphism $\phi_n : \mathbf{Z}[\frac{1}{n}] \to \mathbf{Z}_{n-1}$ of Lemma 2.1, we have that if $x = x_i^k = h(i/n^k)$, then the fixed point in the forward orbit of $i/n^k$ under $\nu_n$ is $\phi_n(i/n^k)$ and the fixed point in the forward orbit of $x$ under $g$ is $h(\phi_n(i/n^k))$. We use $V_j$ with $0 \le j < (n-1)$ to denote the extended orbits of $g$, and we have

(7-6) $$x \in V_j \text{ when } x = h(i/n^k) \text{ and } j = \phi_n(i/n^k).$$

We carry these classes into the endpoints of the partitions and use $V_j P^k$ to denote $V_j \cap VP^k$.

We note that as an element of $\overline{T}_{n,n-1}$ with connected domain, Lemma 2.2 implies that $g$ must cyclically permute the cosets of $p_{n-1}(\Delta_n)$ in $\mathbf{Z}[\frac{1}{n}]$. Since $g$ has fixed points, it fixes one of them, so it fixes them all. Every element of $\mathbf{Z}[\frac{1}{n}]$ has a fixed point of $g$ in its forward orbit, so there is a fixed point of $g$ for each coset. There are as many fixed points as cosets, so there is one fixed point per coset. This discussion applies as well to $\nu_n$ and we get that $h$ induces a well defined permutation on the cosets of $p_{n-1}(\Delta_n)$ in $\mathbf{Z}[\frac{1}{n}]$. The identity coset is preserved by $h$ since we have $h(0) = 0$.

Each interval $I_i^k$ in $P^k$ has a length that we denote by $L(I_i^k)$. These lengths are related. Note that a power of $g$ carries $I_i^s$ affinely onto $I_{i'}^K$ if $s \ge K$ and $i \equiv i' \mod ((n-1)n^K)$. Now (perhaps different) powers of $g$ will carry $I_i^s$ and $I_j^t$ affinely onto the same interval in $P^K$ if $s \ge K$, $t \ge K$ and $i \equiv j \mod ((n-1)n^K)$. Since the affine parts of $g$ have slopes that are integral powers of $n$, this implies that $L(I_i^s)/L(I_j^t)$ is an integral power of $n$.

## 8. Trees.

It is possible to arrange the ingredients of the various partitions $P^k$ into trees. One tree has the endpoints as vertices and the immediate descendents of each $x_i^k$ are all the $x_j^{k+1}$ that map to $x_i^k$ under $g$. Another tree has the intervals as vertices and the immediate descendents of each $I_i^k$ are all the intervals $I_j^{k+1}$ contained in $I_i^k$. The levels of the trees are given by the natural levels of the vertices. The numerical functions that we describe in the following sections are attached to these trees where appropriate. If the reader is helped by viewing the objects as being arranged in trees, then this view is available. The trees are not necessary for the arguments that follow and do not simplify them to any extent. Thus we will not mention these trees in what follows.

## 9. The calculus of break values.

Given an element $x$ of $S_{n-1}$, we associate a measure of the change of slope of $g$ that occurs at $x$. It acts like a derivative, so we give it a notation that resembles $g'$. We define

$$g^b(x) = \log_n \left( \frac{g'_+(x)}{g'_-(x)} \right)$$

where $g'_+(x)$ represents the right derivative of $g$ at $x$ (the slope of $g$ immediately to the right of $x$) and $g'_-(x)$ represents the left derivative of $g$ at $x$. Since $g$ is PL, this is defined for all $x$. We call $g^b(x)$ the break value of $g$ at $x$. We use the logarithm to make the chain rule additive in that if $g_1$ and $g_2$ are two PL functions, then

(9-1) $$(g_1 \circ g_2)^b(x) = g_1^b(g_2(x)) + g_2^b(x).$$

This follows because left and right derivatives satisfy the usual chain rule. We use log base $n$ to make the numbers come out nicer — a fact that we will never use except in examples.

We can think of the break value as the non-linearity concentrated at a point. For a continuously twice differentiable function $f$, the non-linearity $Nf(x)$ of $f$ at $x$ is usually defined to be

$$Nf(x) = (\ln(f'(x)))' = \frac{f''(x)}{f'(x)}.$$

Now

$$\int_a^b Nf = \ln(f')|_a^b$$
$$= \ln(f'(b)) - \ln(f'(a))$$
$$= \ln\left(\frac{f'(b)}{f'(a)}\right).$$

Except for the different log base, this last expression gives the break value at $x$ of a PL function $f$ if we apply the expression to an interval $[a, b]$ with $x \in (a, b)$ and with no other break of $f$ in $[a, b]$. The behavior of the break value and the use that we make of it are similar to the behavior and use of non-linearity.

Some properties of the break value are as follows.

(9-2) $\quad\quad g^b(x) = 0$ if and only if $g$ is affine in a neighborhood of $x$.

$$\text{(9-3)} \qquad (g^n)^b(x) = \sum_{i=o}^{n-1} g^b(g^i(x)).$$

This is just the sum of the break values of $g$ at $n$ consecutive points in an orbit starting with $x$. If $f : S_{n-1} \to S_{n-1}$ is a PL homeomorphism, then

$$\text{(9-4)} \qquad (fgf^{-1})^b(x) = -f^b(f^{-1}(x)) + g^b(f^{-1}(x)) + f^b((gf^{-1})(x))$$

and

$$\text{(9-5)} \qquad (fgf^{-1})^b(f(x)) = -f^b(x) + g^b(x) + f^b(g(x))$$

where (9-4) and (9-5) have been simplified by using the fact that $(f^{-1})^b(f(x)) = -f^b(x)$ which follows from applying (9-1) to $(f^{-1}f)^b(x)$.

Since $g$ is a PL function defined on a circle $S_{n-1}$, the sum of the break values of $g$ once around $S_{n-1}$ must be zero. Since $VP^k$ contains all the break points of $g$ as long as $k \geq K$, we get

$$\text{(9-6)} \qquad \sum_{x \in VP^k} g^b(x) = 0 \text{ whenever } k \geq K.$$

We will use these properties to make two observations, but first we need to understand the break values of $g$ at its fixed points.

From Lemma 6.3, we know that $h$ is a normalizer of $BT_{n,n-1}$. Recall that we denote the groupoid of germs generated by $BT_{n,n-1}$ by $\mathcal{B}_{n,n-1}$. From Lemma 3.3, we know that $h$ induces an isomorphism from the group of germs $(\mathcal{B}_{n,n-1})_{i,i}$ of functions of $BT_{n,n-1}$ with fixed points at the integer $i$ in $S_{n-1}$ to the group of germs $(\mathcal{B}_{n,n-1}))_{h(i),h(i)}$. Now this group of germs is isomorphic to $\mathbf{Z} \times \mathbf{Z}$ and the two natural generators are "slope $n$ to the right of $i$" and "slope $n$ to the left of $i$." Conjugation by an orientation preserving homeomorphism keeps the "behavior to the right" separate from the "behavior to the left" in that the germ to the right of $h(i)$ of the conjugate is determined only by the germ to the right of $i$ of the conjugated function, and similarly for germs to the left. Thus a germ to the right of $i$ that represents "slope $n$ to the right" must be carried to "slope $n$ to the right" or "slope $1/n$ to the right." However, conjugation preserves the property of a fixed point being repelling. Thus "slope $n$ to the right" must be carried to "slope $n$ to the right." Similarly for the slopes to the left. We now note that $\nu_n$ has slope $n$ at all its fixed points — the integer points $i$ in $S_{n-1}$. Thus $g$ has slope $n$ at both sides of all its fixed points $h(i)$. We now can say

$$\text{(9-7)} \qquad g^b(h(i)) = 0, \text{ where the } h(i), \ 0 \leq i < n \text{ are the fixed points of } g.$$

We can now give our first observation. For any $x \in \mathbf{Z}[\frac{1}{n}]$ we have $g^N(x) = h(i)$ is a fixed point of $g$ where $x = h(i/n^N)$ and $N = \lambda(x)$ is the natural level of $x$. Thus every $n \geq N$ has $g^n(x) = h(i)$. Since $g^b(h(i)) = 0$, we get from (9-3) that $(g^n)^b(x)$ is independent of $n$ as long as $n \geq N$. For any $x \notin \mathbf{Z}[\frac{1}{n}]$, the forward orbit of $x$ under $g$ is disjoint from $\mathbf{Z}[\frac{1}{n}]$. Since $\mathbf{Z}[\frac{1}{n}]$ contains all the break points of $g$, we get that $(g^n)^b(x) = 0$ for all $n$ if $x \notin \mathbf{Z}[\frac{1}{n}]$. This allows us to define a function $\Sigma : S_{n-1} \to \mathbf{Z}$ by

$$\text{(9-8)} \qquad \Sigma(x) = \lim_{n \to \infty} (g^n)^b(x).$$

We can refer to it as the total iterated break value at $x$, and it is the sum of all the break values in the forward orbit of $x$ (including $x$) under $g$. We know that the set of breakpoints of $g$ is contained in the set of endpoints of the intervals of $P^K$ where $K$ is the stable level for $g$. Thus all $x \in \mathbf{Z}[\frac{1}{n}]$ with natural level $\lambda(x) > K$ have $g^b(x) = 0$. If $x \in \mathbf{Z}[\frac{1}{n}]$ has $\lambda(x) > K$, then all $(g^i)^b(x) = 0$ for $i < \lambda(x) - K$ and we have

$$\text{(9-9)} \qquad \Sigma(x) = \Sigma(x') \text{ where } x' = g^{\lambda(x)-K}(x) \text{ has } \lambda(x') = K.$$

Since there are only finitely many endpoints of $P^K$ (those $x \in \mathbf{Z}[\frac{1}{n}]$ with $\lambda(x) \leq K$), there are only finitely many different values that $\Sigma(x)$ can achieve. The property given in (9-9) is the reason for calling the level $K$ the stable level.

Note that (9-8) and (9-3) imply that $g^b(x) = \Sigma(x) - \Sigma(g(x))$ for any $x$. This allows us to to say that the existence of the function $\Sigma()$ implies that $g^b$ is a coboundary on an appropriate complex. In smooth settings, this is often enough to reach conclusions about the conjugating homeomorphism $h$. We show in Example 1 of Section 13, that this is not sufficient here.

The second observation is a refinement of the first observation and motivates what is needed beyond the existence of the function $\Sigma()$ to prove that $h$ is PL. We assume that $h$ is PL. We know from Lemma 6.3 that this will imply that $h$ is in $T_{n,n-1}$ so that the break points of $h$ are in $\mathbf{Z}[\frac{1}{n}]$. For $x \in \mathbf{Z}[\frac{1}{n}]$ we have $\nu_n^k(x) = i$ for some integer $i$ with $0 \leq i < n$ and all sufficiently large $k$. Since $\nu_n$ has no break points, we have

$$\text{(9-10)} \qquad (g^k)^b(h(x)) = (h\nu_n^k h^{-1})^b(h(x)) = -h^b(x) + h^b(\nu_n^k(x)) = -h^b(x) + h^b(i)$$

for all sufficiently large $k$. However, $(g^k)^b(h(x)) = \Sigma(h(x))$ for all sufficiently large $k$. Since there are only finitely many break points of $h$, we have for all but finitely many $x \in \mathbf{Z}[\frac{1}{n}]$ that $\Sigma(h(x)) = h^b(i)$ and depends only on the extended orbit of $h(x)$. Thus in the case that $h$ is PL, the observation made above that $\Sigma()$ takes on only finitely many values is further refined. Since we are trying to show that $h$ is PL, we do not need this observation, but it motivates the converse that we prove in the next section.

## 10. Criteria for piecewise linearity.

In this section we develop separate criteria for the piecewise linearity of $h$ and $h^{-1}\nu_n h$. Section 11 has the task of relating these two criteria.

### 10.1. The piecewise linearity of $h$.
We will show that the observation made at the end of the previous section has a converse. We know that if $h$ is PL, then it must have certain properties. If $h$ is not known to be PL, but is known to have these properties, then we will build a PL function $\bar{h}$ that has these properties. Then we will show that it must be equal to $h$ by showing that $\bar{h}$ also conjugates $\nu_n$ to $g$.

From (9-10), we know that if $h$ is PL, if $x \in \mathbf{Z}[\frac{1}{n}]$, and if $i = \phi_n(x)$ is the fixed point of $\nu_n$ in the forward orbit of $x$ under $\nu_n$, then $\Sigma(h(x)) = h^b(i) - h^b(x)$. For a PL function $h$, we know that $h^b(x) = 0$ for almost all $x$ and that will make $\Sigma(h(x))$ take one value $\Sigma_i$ for almost all $x$ with $\phi_n(x) = i$. Since we will have $\Sigma_i = \Sigma(h(x)) = h^b(i) - 0 = h^b(i)$ for those $x$ with $h^b(x) = 0$, we will have $h^b(i) = \Sigma_i$. For those $x$ with $h^b(x) \neq 0$, we will have $h^b(x) = \Sigma_i - \Sigma(h(x))$. This motivates our proof the converse.

Let $h$ be as given in the hypothesis of Theorem 3 and assume, for each $i$ with $0 \leq i < (n-1)$, that $\Sigma(h(x))$ takes on one value $\Sigma_i$ for all but finitely many $x \in \mathbf{Z}[\frac{1}{n}]$ with $\phi_n(x) = i$. For each $x \in \mathbf{Z}[\frac{1}{n}]$, define $b(x) = \Sigma_i - \Sigma(x)$ where $i = \phi_n(x)$. Note that our hypothesis and construction imply that $b(x) = 0$ for all but finitely many $x \in \mathbf{Z}[\frac{1}{n}]$. We want to argue that $h$ is PL with all its break points in $\mathbf{Z}[\frac{1}{n}]$ and that $h^b(x) = b(h(x))$ for all $x \in \mathbf{Z}[\frac{1}{n}]$. We will do this by building an $\bar{h}$ with these properties and then showing that $\bar{h}$ must equal $h$. As in (9-6), we must also have

$$(10\text{-}1) \qquad \sum_{x \in \mathbf{Z}[\frac{1}{n}]} b(h(x)) = \sum_{x \in \mathbf{Z}[\frac{1}{n}]} h^b(x) = 0.$$

in order for this to work. This is a tedious calculation that is based on (9-6) that we now tackle. Note that since $h$ is a one to one correspondence from $\mathbf{Z}[\frac{1}{n}]$ to itself, we can replace $b(h(x))$ with $b(x)$ in (10-1).

We have $b(x) = \Sigma_i - \Sigma(x)$ where $i = \phi_n(x)$. The condition $i = \phi_n(x)$ is equivalent to $h(x) \in V_i$, so we have $b(x) = \Sigma_i - \Sigma(x)$ where $x \in V_i$. We know that for $x \in V_i$, and $\lambda(x) \geq K$, we have $\Sigma(x) = \Sigma_i$. Thus $b(x) = 0$ for all $x$ with $\lambda(x) \geq K$. In spite of the fact that $x$ with $\lambda(x) = K$ contributes nothing to a sum of the values of $b(x)$, it is convenient to include them. We now know that to verify (10-1), it is sufficient to verify

$$\sum_{x \in VP^K} b(x) = \sum_{i=0}^{n-2} \sum_{x \in V_i P^K} b(x) = \sum_{i=0}^{n-2} \sum_{x \in V_i P^K} (\Sigma_i - \Sigma(x)) = 0.$$

We will need to know the sizes of some collections of elements of $\mathbf{Z}[\frac{1}{n}]$. We let $M(k)$ be the number of elements of $VP^k$, we let $M_i(k)$ be the number of elements of $V_iP^k$ and we let $m_i(k)$ be the number of elements $x$ of $V_iP^k$ that have $\lambda(x) = k$. We have

$$M(k) = (n-1)n^k,$$
$$M_i(k) = n^k, \text{ and}$$
$$m_i(k) = n^k - n^{k-1} = (n-1)n^{k-1}.$$

Now
$$\sum_{i=0}^{n-2} \sum_{x \in V_i P^K} (\Sigma_i - \Sigma(x)) = \sum_{i=0}^{n-2} M_i(K)\Sigma_i - \sum_{i=0}^{n-2} \sum_{x \in V_i P^K} \Sigma(x)$$
$$= n^K \sum_{i=0}^{n-2} \Sigma_i - \sum_{x \in VP^K} \Sigma(x)$$

so it suffices to prove
$$\sum_{x \in VP^K} \Sigma(x) = n^K \sum_{i=0}^{n-2} \Sigma_i.$$

We have

(10-2) $$\sum_{x \in VP^K} \Sigma(x) = \sum_{\lambda(x)=K} \Sigma(x) + \sum_{x \in VP^{K-1}} \Sigma(x).$$

But $x \in V_i P^K$ and $\lambda(x) = K$ imply $\Sigma(x) = \Sigma_i$ so we have

(10-3) $$\sum_{\lambda(x)=K} \Sigma(x) = m_i(K) \sum_{i=0}^{n-2} \Sigma_i = (n-1)n^{K-1} \sum_{i=0}^{n-2} \Sigma_i.$$

Thus we have
$$\sum_{x \in VP^K} \Sigma(x) = (n-1)n^{K-1} \sum_{i=0}^{n-2} \Sigma_i + \sum_{x \in VP^{K-1}} \Sigma(x)$$

and we are done if we show
$$\sum_{x \in VP^{K-1}} \Sigma(x) = n^{K-1} \sum_{i=0}^{n-2} \Sigma_i.$$

By (10-3), this is equivalent to
$$(n-1) \sum_{x \in VP^{K-1}} \Sigma(x) = \sum_{\lambda(x)=K} \Sigma(x)$$

and from (10-2) this is equivalent to the following which will be more convenient to prove:

$$(10\text{-}4) \qquad n \sum_{x \in VP^{K-1}} \Sigma(x) = \sum_{x \in VP^K} \Sigma(x).$$

We now expand the expression $\Sigma(x)$. From (9-8) and (9-3), this is a sum of values of $g^b(g^i(x))$ as $i$ runs over the non-negative integers. The various $g^i(x)$ have natural levels given by $\lambda(g^i(x)) = \lambda(x) - i$ as long as the result is not negative. If we sum all $\Sigma(x)$ over all $x$ with $\lambda(x)$ fixed at some $k$, then we will get contributions of the form $g^b(y)$ from elements $y$ with $\lambda(y) \leq \lambda(x)$. Since there are $n$ elements of level $k$ that map to a given element of level $k-1$, $n^2$ elements of level $k$ that map to a given element of level $k-2$ and so forth, we can assess the multiplicity of contribution of elements of various levels. The pattern breaks down at level 0 however. This is because while a fixed point has $n$ preimages, one of them is the fixed point itself. Thus there are only $(n-1)$ elements of level 1 that map to a given fixed point. We will ignore this because by (9-7), the fixed points do not contribute to the expressions $\Sigma(x)$. Also for this reason, it makes no difference whether our sums include level zero. We attack (10-4) by writing down a sum of $\Sigma(x)$ over a single level, and then summing over levels.

We have

$$\sum_{\lambda(x)=k} \Sigma(x) = \sum_{\lambda(x) \leq k} n^{k-\lambda(x)} g^b(x).$$

Summing over various levels, we get

$$\sum_{\lambda(x) \leq j} \Sigma(x) \; = \; \sum_{k \leq j} \sum_{\lambda(x) \leq k} n^{k-\lambda(x)} g^b(x) \; = \; \sum_{k \leq j} \sum_{i \leq k} \sum_{\lambda(x)=i} n^{k-i} g^b(x)$$

$$= \; \sum_{i=0}^{j} \sum_{\lambda(x)=i} \sum_{k=i}^{j} n^{k-i} g^b(x) \; = \; \sum_{i=0}^{j} \sum_{\lambda(x)=i} \sum_{k=0}^{j-i} n^{k} g^b(x)$$

$$= \; \sum_{i=0}^{j} \sum_{\lambda(x)=i} \left( \frac{n^{j-i+1} - 1}{n-1} \right) g^b(x).$$

Multiplying through by $(n-1)$ gives

$$(n-1) \sum_{\lambda(x) \leq j} \Sigma(x) = \sum_{i=0}^{j} \sum_{\lambda(x)=i} \left( n^{j-i+1} - 1 \right) g^b(x), \quad \text{and}$$

$$(n-1) \sum_{\lambda(x) \leq j-1} \Sigma(x) = \sum_{i=0}^{j-1} \sum_{\lambda(x)=i} \left( n^{j-i} - 1 \right) g^b(x)$$

$$= \sum_{i=0}^{j} \sum_{\lambda(x)=i} \left( n^{j-i} - 1 \right) g^b(x)$$

because $n^{j-j} = 1$.

Breaking these sums into two parts each gives

$$(10\text{-}5) \qquad (n-1) \sum_{\lambda(x) \leq j} \Sigma(x) = \sum_{i=0}^{j} \sum_{\lambda(x)=i} \left(n^{j-i+1}\right) g^b(x) - \sum_{x \in VP^j} g^b(x),$$

and

$$(n-1) \sum_{\lambda(x) \leq j-1} \Sigma(x) = \sum_{i=0}^{j} \sum_{\lambda(x)=i} \left(n^{j-i}\right) g^b(x) - \sum_{x \in VP^j} g^b(x).$$

But from (9-6), we know that

$$\sum_{x \in VP^j} g^b(x) = 0$$

whenever $j \geq K$. Thus we set $j = K$, we factor an $n$ from the right side of (10-5), and we obtain (10-4). This proves (10-1) and we now continue with the construction of $\bar{h}$.

Let $m$ be a real number greater than 0. We build a function $\bar{h}_m$ on $[0, n-1]$ as follows. It essentially "integrates" the information given by the the function $b(x)$. Except for $x = 0$, we draw horizontal lines $H_x$ in the plane through the points $(0, x)$ on the $y$-axis that have $b(x) \neq 0$. We draw a graph of a function from $[0, \infty)$ to itself to satisfy the following: The graph is piecewise linear, it passes through $(0,0)$, it has slope $m$ to the right of 0, it breaks only when it intersects one of the horizontal lines $H_x$, and the break value that it has there is the value $b(x)$. Note that all slopes of this function will be positive. That such a function can be built, and the uniqueness of the function given $m$ is clear. We let $\bar{h}_m$ be the function with this graph. The $x$ coordinate at which $\bar{h}_m$ achieves the value $(n-1)$ is a continuous monotone function of $m$, and in fact it is easily seen to depend linearly on $m^{-1}$. Thus there is a unique $m$ for which $\bar{h}_m(n-1) = (n-1)$. We let $\bar{h} = \bar{h}_m$ for this $m$. We view $\bar{h}$ as a homeomorphism from $S_{n-1}$ to itself. This will make the sum of the break values of $\bar{h}$ over all of $S_{n-1}$ equal to zero. By construction, we have $\bar{h}^b(\bar{h}^{-1}(x)) = b(x)$ at all $x$ except possibly $x = 0$. Because of (10-1), we will have that the break value of $\bar{h}$ at $0 = \bar{h}(0)$ will be $b(0)$.

If we let $\bar{g} = \bar{h} \nu_n \bar{h}^{-1}$, then we can define $\overline{\Sigma}(x)$ as the sum of the breaks in the forward orbit of $x$ under $\bar{g}$ in analogy with $\Sigma(x)$. From (9-10) applied to $\bar{h}$, we get

$$\overline{\Sigma}(\bar{h}(x)) = \bar{h}^b(i) - \bar{h}^b(x)$$

or equivalently

$$\overline{\Sigma}(x) = \bar{h}^b(\bar{h}^{-1}(i)) - \bar{h}^b(\bar{h}^{-1}(x))$$
$$= b(i) - b(x)$$
$$= \left(\Sigma_i - \Sigma(i)\right) - \left(\Sigma_i - \Sigma(x)\right)$$
$$= \Sigma(x)$$

since $\Sigma(i) = 0$. From the definitions of $\overline{\Sigma}()$ and $\Sigma()$, we know that we can recover $\bar{g}^b(x)$ as $\overline{\Sigma}(x) - \overline{\Sigma}(\bar{g}(x))$ and $g^b(x)$ as $\Sigma(x) - \Sigma(g(x))$. From this, we conclude that $\bar{g}^b$ and $g^b$ are identical functions. Further, we know that $\bar{g}(0) = g(0) = 0$ since $\nu_n$, $\bar{h}$ and $h$ all fix 0. Just as the properties of $\bar{h}$ made $\bar{h}$ unique, we can argue that the properties listed make $g$ unique. Thus $\bar{g} = g$ and $\bar{h}$ conjugates $\nu_n$ to $g$. Since $\nu_n$ is an expanding map, the homeomorphism that conjugates $\nu_n$ to $g$ is unique. One way to se this is to note that the conjugator must take the endpoints of the standard Markov partition of level $k$ for $\nu_n$ to the elements of $VP^k$ in order. This forces the value of the conjugator on the elements of $\mathbf{Z}[\frac{1}{n}]$, a dense set in $S_{n-1}$. If we were arguing existence of a conjugator, we would have to argue the denseness of $\bigcup_k \{VP^k\}$. However, a conjugator (namely $h$) exists by hypothesis.

## 10.2. The piecewise linearity of $h^{-1}\nu_n h$.

This title is somewhat misleading. Our aim is to show that $h$ is PL. However, from Lemma 6.3, we know that $h$ normalizes $BT_{n,n-1}$, so not only must $h^{-1}\nu_n h$ be PL, but so must $h^{-1}fh$ for any $f \in BT_{n,n-1}$. Our approach will be to show that if $h$ is not PL (violates the criterion in 10.1), then $h$ must fail to normalize $BT_{n,n-1}$. Being able to test any $f \in BT_{n,n-1}$ gives more flexibility than just testing $h^{-1}\nu_n h$. Thus we develop a criterion for the piecewise linearity of any $h^{-1}fh$ for any $f \in BT_{n,n-1}$.

We can show that a function is not PL by showing that a sequence of chords to its graph changes slope infinitely many times in a bounded interval. This allows us to reduce our calculations of the behavior to a countable set of points. We have enough information about $h$ and $h^{-1}$ from the Markov partitions for $\nu_n$ and $g$ to calculate the behavior of $h^{-1}fh$ on a countable set of points, and this gives us enough machinery to derive our criterion.

Let $P^s$ and $P^t$ be two Markov partitions obtained from standard Markov partitions of levels $s$ and $t$ for $\nu_n$. We know that either $s = t$ or one of $P^s$ or $P^t$ is derived from the other (it doesn't matter which). Assume that $(x_1 < x_2 < \cdots < x_k)$ is a sequence in $VP^s$ so that all points $[x_1, x_k] \cap VP^s$ are in the sequence. Assume that $(x'_1 < x'_2 < \cdots < x'_k)$ is a sequence in $VP^t$ containing all $[x'_1, x'_k] \cap VP^t$. Assume that $f \in BT_{n,n-1}$ takes $x_i$ to $x'_i$ for $1 \leq i < k$ and that it fails to take $x_k$ to $x'_k$. We have that the points $h^{-1}(x_i)$ are equally spaced (they are consecutive endpoints in the standard Markov partition of level $s$ for $\nu_n$) as are the points $h^{-1}(x'_i)$. Thus the graph of $h^{-1}fh$ has chords of the same slope spanning all the intervals $[h^{-1}(x_i), h^{-1}(x_{i+1})]$ for $1 \leq i \leq (k-2)$ but the chord spanning $[h^{-1}(x_{k-1}), h^{-1}(x_k)]$ has a different slope. Thus $h^{-1}fh$ is not affine on $[h^{-1}(x_1), h^{-1}(x_k)]$. In order for this to be a useful observation, we need to show that we can arrange $f$ to take two consecutive endpoints in a Markov partition for $g$ to two consecutive endpoints in a (perhaps different) Markov partition for $g$, and we need to understand when $f$ fails to carry a string of consecutive endpoints to a string of consecutive endpoints.

If $f$ takes two consecutive endpoints from one partition to two consecutive endpoints from another partition, then $f$ is carrying some interval of one partiton to an interval of the other. We achieve this by starting with the fact that $f$ carries an endpoint of one partition to an endpoint of the other and then take derived partitions until we get what we want.

Let $f$ take some $a = x_i^s$ to $b = x_j^t$. Since $f$ preserves cosets of $p_{n-1}(\Delta_n)$, we know that the forward orbits of $a$ and $b$ include the same fixed point of $g$. Specifically, $g^{\lambda(a)}(a) = g^{\lambda(b)}(b) = z$ is a fixed point of $g$ and $\lambda(a)$ and $\lambda(b)$ are the smallest powers of $g$ that achieve this. Assume temporarily that $s = \lambda(a) + K$ and $t = \lambda(b) + K$. Then $g^{\lambda(a)}(I_i^s) = g^{\lambda(b)}(I_j^t) = I_{i'}^K$ where $i \equiv j \equiv i' \mod (n-1)n^K$. We know that for these values of $s$ and $t$, that $L(I_j^t)/L(I_i^s)$ is an integral power of $n$.

If $f$ has this ratio as its slope to the right of $a$, then $f$ carries $I_i^s$ affinely onto $I_j^t$ and we have that $f$ takes two consecutive endpoints of $P^s$ to two consecutive endpoints of $P^t$. For future reference, we call this ratio (with $s = \lambda(a) + K$ and $t = \lambda(b) + K$) the natural slope of the pair $(a, b)$.

If $f$ has slope to the right of $a$ that is different from the natural slope of $(a, b)$, then it differs from it by a factor that is an integral power of $n$. Assume first that the slope of $f$ is too large so that $L(I_i^s)$ is too large. We note that $I_{ni}^{s+1}$ has $a = x_i^s = x_{ni}^{s+1}$ as left endpoint and since $g^{\lambda(a)}(a)$ is a fixed point of $g$, we have $g^{\lambda(a)+1}(I_{ni}^{s+1}) = I_{i'}^K$ in $P^K$ and $g^{\lambda(a)}(I_{ni}^{s+1}) = I_{ni'}^{K+1}$ in $P^{K+1}$ where $I_{i'}^K$ and $I_{ni'}^{K+1}$ have the same left endpoint $z$. Because $z$ is a fixed point of $g$, we have that $g(I_{ni'}^{K+1}) = (I_{i'}^K)$. We know that $g$ has slope $n$ near its fixed points, so $L(I_{i'}^K)/L(I_{ni'}^{K+1}) = n$. The slope of $g^{\lambda(a)}$ on $I_{ni}^{s+1}$ is the same as it has on $I_i^s$ so we can conlude that $L(I_i^s)/L(I_{ni}^{s+1}) = n$ and that the ratio $L(I_j^t)/L(I_{ni}^{s+1})$ differs from the natural slope of $(a, b)$ by a factor of $1/n$.

By taking further derivations of $P^s$, we can achieve an adjustment by any integral multiple of $1/n$ we wish. If the slope of $f$ is too small, then we take derivations of $P^t$ instead of $P^s$. In this way we can get $f$ to take two consecutive points of one partition to two consecutive points of another. Further, we can do this with any $f$ in $BT_{n,n-1}$ and have as the two first points any $a$ and $b$ in $\mathbf{Z}[\frac{1}{n}]$ where $f(a) = b$.

Now assume that an element $f$ of $BT_{n,n-1}$ carries the sequence of consecutive endpoints $(a = x_i^s, x_{i+1}^s, \ldots, x_{i+k}^s = b)$ in a partition $P^s$ for $g$ to consecutive endpoints $(c = x_j^t, x_{j+1}^t, \ldots, x_{j+k}^t = d)$ in a partition $P^t$ for $g$. Assume further that $f$ is affine on $[a, b]$. If both $P^s$ and $P^t$ have their derivations taken, then $f$ carries the new endpoints in $[a, b]$ in $P^{s+1}$ onto the new endpoints in $[c, d]$ in $P^{t+1}$. Thus we can derive domain and range partitions the same number of times and keep the assumed hypotheses on the original interval $[a, b]$. Thus we can assume that $s \geq \lambda(a) + K$ and $t \geq \lambda(c) + K$. As observed before, $a$ and $c$ have a common fixed point $z$ of $g$ in their forward orbits. Our assumptions on $s$ and $t$ guarantee that for this fixed point $z$ of $g$, we have $g^{s-K}(a) = g^{t-K}(c) = z = x_{i'}^K$ where $i' \equiv i \equiv j \mod (n-1)n^K$.

We now assume that $f$ is affine on $[a, x_{i+k+1}^s]$. We wish to see if $f$ takes $x_{i+k+1}^s$ to $x_{j+k+1}^t$. This is equivalent to having $f$ carry $I_{i+k}^s$ onto $I_{j+k}^t$. This in turn is equivalent to having the lengths of $I_{j+k}^t$ and $I_{i+k}^s$ have the same ratio as the lengths of $I_{j+k-1}^t$ and $I_{i+k-1}^s$. The intervals $I_{i+k-1}^s$ and $I_{i+k}^s$ are carried affinely by $g^{s-K}$ onto $I_{i'+k-1}^K$ and $I_{i'+k}^K$ with slopes $m_1$ and $m_2$ respectively. The intervals $I_{j+k-1}^t$ and $I_{j+k}^t$ are carried affinely by $g^{t-K}$ onto $I_{i'+k-1}^K$ and $I_{i'+k}^K$ with slopes $m_3$ and $m_4$ respectively. Our chain of equivalent statements is equivalent to $m_3/m_1 = m_4/m_2$ which is equivalent to $m_2/m_1 = m_4/m_3$. Thus we have the equality of break values

$$(g^{s-K})^b(x_{i+k}^s) = (g^{t-K})^b(x_{j+k}^t). \tag{10-6}$$

But

$$g^{s-K}(x_{i+k}^s) = g^{t-K}(x_{j+k}^t) = x_{i'+k}^K.$$

So (10-6) is equivalent to

$$\Sigma(x_{i+k}^s) = (g^s)^b(x_{i+k}^s) = (g^{s-K})^b(x_{i+k}^s) + (g^K)^b(x_{i'+k}^K)$$
$$= (g^{t-K})^b(x_{j+k}^t) + (g^K)^b(x_{i'+k}^K) = (g^t)^b(x_{j+k}^t) = \Sigma(x_{j+k}^t).$$

Thus extending the fact that an affine $f$ carries a sequence of consecutive endpoints to a sequence of consecutive endpoints is equivalent to the equality of the values of $\Sigma(x)$ on the corresponding endpoints (except the last endpoint in each sequence). Specifically, given that $f$ is affine on $[a, b]$ as given above, and takes the pair $(x_i^s, x_{i+1}^s)$ to the pair $(x_j^t, x_{j+1}^t)$ with $s \geq \lambda(a) + K$ and $t \geq \lambda(c) + K$, then it takes $(x_i^s, x_{i+1}^s, \ldots, x_{i+k}^s)$ to $(x_j^t, x_{j+1}^t, \ldots, x_{j+k}^t)$ if and only if all of $\Sigma(x_{i+q}^s) = \Sigma(x_{j+q}^t)$ for all $q$ with $1 \leq q < k$.

## 11. Relating the criteria.

In this section, we restrict $n$ to 2. Theorem 3 follows from the next proposition, from the criterion for the piecewise linearity of $h$, and the fact that an $f$ can be built as described in the proposition by Lemma 2.3. Recall that the $V_i$, $0 \leq i < (n-1)$, are the extended oribts of $g$ in $\mathbf{Z}[\frac{1}{n}]$. When $n = 2$, there is only one such set $V_0 = \mathbf{Z}[\frac{1}{2}]$. The criterion in Section 10.1 that $h$ is PL if the function $\Sigma()$ is almost constant on each $V_i$ reduces to requiring that $\Sigma()$ be constant off a finite set in $\mathbf{Z}[\frac{1}{2}]$.

PROPOSITION 11.1. *Assume that $\Sigma()$ is not constant on the complement of a finite set in $\mathbf{Z}[\frac{1}{2}]$. Let $a$ and $c$ be in $\mathbf{Z}[\frac{1}{2}]$ with $s = \lambda(a) \geq K$ and $t = \lambda(c) \geq K$. Let $a = x_j^s$ and $b = x_{j+1}^s$. If $f \in BT_{2,1}$ takes $a$ to $c$, is affine on $[a,b]$ and does not have the natural slope of the pair $(a,c)$ on the interval $[a,b]$, then $h^{-1}fh$ is not PL on $[h^{-1}(a), h^{-1}(b)]$.*

The proposition follows immediately by induction from the next lemma.

LEMMA 11.2 (INDUCTIVE STEP). *If $f$, $\Sigma()$, $[a,b]$ and $c$ satisfy the hypotheses of the proposition, then there is an interval $[d,e] \subseteq (a,b]$ so that $f$ and $[d,e]$ satisfy the hypotheses of the proposition and so that $d$ is a break point of $h^{-1}fh$.*

PROOF: The bulk of the proof is in finding the point $d$ at which $h^{-1}fh$ breaks. After that it is a straightforward matter to verify that $d$ has the remaining properties claimed.

**11.1. Locating $d$.** We have that $f$ takes $a$ to $c$. We have $a = x_j^s$ where $s = \lambda(a)$. With $t = \lambda(c)$, we can say $c = x_k^t$. From Section 10.2, we have that the natural slope of $(a,c)$ is gotten by regarding $a = x_{2^K j}^{s+K}$ and $c = x_{2^K k}^{t+K}$ and that the natural slope is equal to $L(I_{2^K k}^{t+K})/L(I_{2^K j}^{s+K})$. Given a positive integer $p$, we also know that

(11-1) $$L(I_{2^{K+p} j}^{s+K+p}) = 2^{-p} L(I_{2^K j}^{s+K}) \text{ and } L(I_{2^{K+p} k}^{t+K+p}) = 2^{-p} L(I_{2^K k}^{t+K}).$$

Since the slope $m$ of $f$ on $[a,b]$ is not the natural slope of $(a,c)$, we know that either $m = L(I_{2^{K+p} k}^{t+K+p})/L(I_{2^K j}^{s+K})$ or $m = L(I_{2^K k}^{t+K})/L(I_{2^{K+p} j}^{s+K+p})$ for some positive integer $p$. From this it follows that either $f(I_{2^k j}^{s+K}) = I_{2^{K+p} k}^{t+K+p}$ or $f(I_{2^{K+p} j}^{s+K+p}) = I_{2^K k}^{t+K}$ and $f$ takes two consecutive endpoints of $P^{s+K}$ or $P^{s+K+p}$ to two consecutive endpoints of $P^{t+K+p}$ or $P^{t+K}$ respectively.

To verify the criterion in Section 10.2, we must compare $\Sigma()$ on the endpoints of $P^{s+K}$ or $P^{s+K+p}$ to corresponding endpoints of $P^{t+K+p}$ or $P^{t+K}$ respectively. Since the direction of $f$ is irrelevant to the comparison, we will assume that $m = L(I_{2^{K+p} k}^{t+K+p})/L(I_{2^K j}^{s+K})$. Later we will say a few words about the other possibility.

When we locate $d$, we will want to claim that $d$ is a break point of $h^{-1}fh$. We note that our criterion in Section 10.2 only detects breaks in the slopes of chords to the graph of $h^{-1}fh$. A break in the slopes of chords implies that a break in the graph of $h^{-1}fh$ takes place somewhere in the domain covered by the two chords with differing slopes, but it does not precisely locate the break of $h^{-1}fh$. To make the location more precise, we will investigate the slopes of chords on increasingly fine intervals. The next paragraph shows how we can do this.

We are assuming that $f(I^{s+K}_{2^K j}) = I^{t+K+p}_{2^{K+p} k}$. This means that $f(x^{s+K}_{2^K j+i}) = x^{t+K+p}_{2^{K+p} k+i}$ for $i = 0$ and $i = 1$. We are to see if $\Sigma(x^{s+K}_{2^K j+i}) = \Sigma(x^{t+K+p}_{2^{K+p} k+i})$ for various values of $i$. However, (11-1) and our assumption that $f(I^{s+K}_{2^K j}) = I^{t+K+p}_{2^{K+p} k}$ imply that

(11-2) $$f(I^{s+K+q}_{2^{K+q} j}) = I^{t+K+p+q}_{2^{K+p+q} k}$$

for any non-negative integer $q$. Thus we have $f(x^{s+K+q}_{2^{K+q} j+i}) = x^{t+K+p+q}_{2^{K+p+q} k+i}$ for $i = 0$ and $i = 1$ and we can also see if $\Sigma(x^{s+K+q}_{2^{K+q} j+i}) = \Sigma(x^{t+K+p+q}_{2^{K+p+q} k+i})$ for various values of $i$. In other words, we can not only compare values of $\Sigma()$ on sequences in $VP^{s+K}$ and $VP^{t+K+p}$, but also compare values of $\Sigma()$ on sequences in $VP^{s+K+q}$ and $VP^{t+K+p+q}$.

Note that $VP^{s+K} \subseteq VP^{s+K+q}$ and $VP^{t+K+p} \subseteq VP^{t+K+p+q}$. What we will show will be that a discrepancy in the values of $\Sigma()$ will occur at some point in $VP^{s+K} \cap (a, b)$, but that no discrepancy occurs in $(a, b)$ at a point in $VP^{s+K+q} - VP^{s+K}$. Since $VP^{s+K}$ is finite, any discrepancy in $(a, b)$ in any of the $VP^{s+K+q}$ occurs in a fixed finite set of points. We will let $d$ be the leftmost point in $(a, b)$ at which a discrepancy occurs.

Since no discrepancy occurs to the left of $d$ in any $VP^{s+K+q}$, we have that the slopes of chords of $h^{-1}fh$ are the same on a dense set of points in $[a, d]$ so $h^{-1}fh$ is affine on $[h^{-1}(a), h^{-1}(d)]$. If $d'$ is the point immediately to the right of $d$ in $VP^{s+K}$, then we know that $h^{-1}fh$ is not affine on $[h^{-1}(a), h^{-1}(d')]$. However, if $h^{-1}fh$ is affine on an interval in $[h^{-1}(a), h^{-1}(d')]$ that properly contains $[h^{-1}(a), h^{-1}(d)]$, then there would have to be a discrepancy of values of $\Sigma()$ at some point $d''$ in $(d, d')$ in some $VP^{s+K+q}$. But since $d$ and $d'$ are consecutive points in $VP^{s+K}$, we would have to have $d'' \in VP^{s+K+q} - VP^{s+K}$. Since this will be shown not to happen, we must have that $h^{-1}fh$ breaks at $d$.

We are to see if the equalities $\Sigma(x^{s+K+q}_{2^{K+q} j+i}) = \Sigma(x^{t+K+p+q}_{2^{K+p+q} k+i})$ hold. Since the notation is getting out of hand, and since we will need to understand the effect of powers of $g$ on the indexes, we shift our point of view. We make use of (7-3) which gives

$$x^{s+K+q}_{2^{K+q} j+i} = h((2^{K+q} j + i)/2^{s+K+q})$$

and

$$x^{t+K+p+q}_{2^{K+p+q} k+i} = h((2^{K+p+q} k + i)/2^{t+K+p+q}).$$

This will allow us to refer to the binary expansions of $(2^{K+q} j + i)/2^{s+K+q}$ and $(2^{K+p+q} k + i)/2^{t+K+p+q}$.

Let $x = h(a/2^b)$ where $a$ and $b$ are positive integers. We know from the remarks following (7-3) that $\lambda(x)$ is the position to the right of the binary point of the last 1 in the binary expansion of $a/2^b$. Assume that $\lambda(x) \geq K$. From (9-9), we know that $\Sigma(x) = \Sigma(g^{\lambda(x)-K}(x))$. Now

(11-3) $$g^{\lambda(x)-K}(x) = h\nu_2^{\lambda(x)-K}h^{-1}(x) = h\nu_2^{\lambda(x)-K}(a/2^b) = h\Omega_K(a/2^b),$$

where $\Omega_K(a/2^b)$ is defined to be $\nu_2^{\lambda(x)-K}(a/2^b)$. If we refer to the last 1 in the binary expansion of $a/2^b$ (which we know occurs in position $\lambda(x)$) as the "last" bit in the expansion of $a/2^b$, then the expansion of $\Omega_K(a/2^b)$ consists of a binary point followed by the "last" $K$ bits of the expansion of $a/2^b$. Thus $\Sigma(x)$ depends on the "last" $K$ bits of the binary expansion of $h^{-1}(x) = a/2^b$.

It will be convenient to introduce another notational device. For an integer $i$, we let $z(i)$ be the largest integer $z$ for which $2^z \mid i$ holds. This is the number of trailing zeros in the binary expansion of $i$.

We consider a pair of points $x$ and $y$ at which values of $\Sigma()$ are to be compared. We let $x = h((2^{K+q}j+i)/2^{s+K+q})$. Since $s = \lambda(a) = \lambda(x_j^s)$, we know that the last 1 in the binary expansion of $j/2^s$ is in position $s$ to the right of the binary point. We are only interested in points in $(a,b) = (x_j^s, x_{j+1}^s)$ and $x_{j+1}^s = h((j+1)/2^s)$. Now $(2^{K+q}j + i)/2^{s+K+q} = (j/2^s) + (i/2^{s+K+q})$. Thus we are only interested in values of $i$ with $1 \leq i < 2^{K+q}$, and $\lambda(x) = s + K + q - z(i)$. The point $y$ corresponding to $x$ is given by

$$y = h((2^{K+p+q}k + i)/2^{t+K+p+q}) = h((k/2^t) + (i/2^{t+K+p+q})),$$

and we have $\lambda(y) = t + K + p + q - z(i)$.

As $i$ runs through its range of values, then $x$ and $y$ run through consecutive points in $VP^{s+K+q}$ and $VP^{t+K+p+q}$ that are in $(a,b)$ and $(f(a), f(b))$ respectively. However, we can break this into levels, and observe that as $i$ runs through those values in its range that have $z(i)$ equal to some constant $\zeta$, then $x$ runs, in order, through the points in $VP^{s+K+q} \cap (a,b)$ that are at the constant natural level $s + K + q - \zeta$, and similarly $y$ runs, in order, through those points in $VP^{t+K+p+q} \cap (f(a), f(b))$ at the natural level $t + K + p + q - \zeta$. The crucial point is that as $i$ runs through its values with $z(i) = \zeta$, the points $g(x)^{\lambda(x)-K}$ are running, in order, through a set of consecutive point on level $K$, and the points $g^{\lambda(x)-K}(y)$ are running through a perhaps different set of consecutive points also on level $K$. We will exploit the fact that values of $\Sigma()$ will be compared on different sequences of points on level $K$. We must now compute which sequences will be compared.

Let $\zeta$ be fixed and let $x = h((j/2^s) + (i/2^{s+K+q}))$ where $z(i) = \zeta$. We have $\lambda(x) = s + K + q - \zeta$. From (9-9) and (11-3), we know that

$$\Sigma(x) = \Sigma\Big(h\Omega_K((j/2^s) + (i/2^{s+K+q}))\Big).$$

Now the last 1 in the expansion of $j/2^s$ is at position $s$ and the last 1 in the expansion of $i/2^{s+K+q}$ is at position $s+K+q-\zeta$. Thus if $\lambda(x) = s+K+q-\zeta \geq s+K$, then $\Omega_K((j/2^s)+(i/2^{s+K+q}))$ depends only on $i$ and so does $\Sigma(x)$. Further, the definition of $\Omega_K$ makes it clear that

$$\Omega_K((j/2^s)+(i/2^{s+K+q})) = \Omega_K(i/2^{s+K+q}) = \Omega_K(i/2^K) = i/2^K \mod (1),$$

and we get that $\Sigma(x) = \Sigma\big(h\Omega(i/2^K)\big)$ whenever $\lambda(x) \geq s+K$. Similarly, for $y = h((k/2^t)+(i/2^{t+K+p+q}))$, we have $\Sigma(y) = \Sigma\big(h\Omega(i/2^K)\big)$ whenever $\lambda(y) \geq t+K$. This establishes that no discrepancy in the values of $\Sigma$ will be discovered when we are considering points $x$ from $VP^{s+K+q} - VP^{s+K}$. (In fact this is true for $VP^{s+K+q} - VP^{s+K-1}$, but there is no need to dwell on this.)

From now on we take $q$ to be 0. We are looking at points $x = h((j/2^s)+(i/2^{s+K}))$ and $y = h((k/2^t)+(i/2^{t+K+p}))$ where $1 \leq i < 2^K$. The restrictions on $i$ dictate that the binary expansions of $i/2^{s+K}$ and $i/2^{t+K+p}$ have no ones outside the "last" $K$ bits of the expansion. We have $\lambda(x) = s+K-z(i)$ and $\lambda(y) = t+K+p-z(i)$. We now look at all those $i$ with $z(i)$ a fixed value $\zeta \geq 0$.

For these $i$, we have have

$$\Sigma(x) = \Sigma\big(h\Omega_K((j/2^s)+(i/2^{s+K}))\big)$$

where $\Omega_K((j/2^s)+(i/2^{s+K}))$ expands as the "last" $\zeta$ bits of $j/2^s$ followed by the $K-\zeta$ bits of $i/2^{s+K}$ that end at the last 1 in the expansion of $i/2^{s+K}$. For a similar $i$, we have

$$\Sigma(y) = \Sigma\big(h\Omega_K((k/2^t)+i/2^{t+K+p}))\big)$$

where $\Omega_K((k/2^t)+i/2^{t+K+p}))$ expands as the "last" $\zeta - p$ bits of $k/2^t$ followed by $p$ zeros followed by the same $K - \zeta$ bits of $i/2^{s+K}$ that end at the last 1 in the expansion of $i/2^{s+K}$ that are the last $K - \zeta$ bits of $\Omega_K((j/2^s)+(i/2^{s+K}))$.

Thus we see that the sequence of values $\Sigma(x)$ that we get for those $i$ with $z(i) = \zeta$ are those $\Sigma\big(h(w)\big)$ where $w$ runs over the $2^{K-\zeta}$ patterns of $K$ bits that start with the last $\zeta$ bits of $j/2^s$. This means that all the patterns $w$ have a common prefix of $\zeta$ bits and that the bit in position $\zeta$ (the last bit of $j/2^s$) is a 1. The sequence of values $\Sigma(y)$ that we get are those $\Sigma\big(h(v)\big)$ where $v$ runs over the $2^{K-\zeta}$ patterns of $K$ bits that start with a prefix of $\zeta$ bits and where the bit in position $\zeta$ is a 0.

Note that our use of (9-9) means that all values of $\Sigma()$ are reduced to values of $\Sigma(u)$ with $\lambda(u) = K$. Since there are $2^{K-1}$ different $u$ with $\lambda(u) = K$ (see Section 10.1), we know that the full sequence of values of $\Sigma()$ that are available is a sequence with $2^{K-1}$ elements. We call a string of $2^{K-\zeta}$ such values given by all $\Sigma(h(w))$ with $w$ a $K$ bit expression with a fixed $\zeta$ bit prefix a $\zeta$-string of the values of $\Sigma()$. The full sequence of $2^{K-1}$ possible $\Sigma()$ values is divided into $2^\zeta$ different $\zeta$-strings.

We call a $\zeta$-string an even $\zeta$-string if the prefix of $\zeta$ bits that determines the $\zeta$-string ends in a 0. Similarly, we call a $\zeta$-string an odd $\zeta$-string if the prefix of $\zeta$ bits that determines the $\zeta$-string ends in a 1. Note that remarks made above say that we are always comparing an even $\zeta$-string with an odd $\zeta$-string.

Note that $\zeta$ can be fixed at any value from 0 through $K-1$. At $\zeta = 0$, we are comparing all of the sequence of $2^{K-1}$ values of $\Sigma()$ with itself. This results in no discrepancy and reflects that parenthetical remark above that no discrepancy will be discovered for an $x$ outside of $VP^{s+K-1}$. For $\zeta = 1$ we compare the unique odd 1-string with the unique even 1-string. In simpler words, we compare the first half of the sequence of $\Sigma()$ values to the second half of the sequence. For $\zeta = 2$, we compare some odd quarter to some even quarter. The pattern continues.

We now have a very elementary induction. If the first half is the same as the second half (as sequences), then all odd quarters are equal to each other (as sequences) and all even quarters are equal to each other (as sequences). If now one odd quarter is shown to equal one even quarter (as sequences), then all four quarters are equal (as sequences). We now note by induction that the sequence of $2^{K-1}$ values of $\Sigma()$ must be constant if all comparisons of some even $\zeta$-string with some odd $\zeta$-string show no discrepancies for all $\zeta$ with $1 \leq \zeta \leq (K-1)$. Since we assume that $\Sigma()$ is not constant over those $u$ with $\lambda(u) = K$, we must have a discrepancy at some level.

This locates $d$ under the assumption that $m = L(I_{2^K+p_k}^{t+K+p})/L(I_{2^K j}^{s+K})$.

If we had assumed that $m = L(I_{2^K k}^{t+K})/L(I_{2^K+p_j}^{s+K+p})$, then an identical argument would locate a discrepancy at some $x = h((j/2^s) + (i/2^{s+K+p}))$ with $1 \leq i < 2^K$. This will be a point that is well within the interval $(x_j^s, x_{j+1}^s)$.

**11.2. Verifying properties.** We will work with the point $d$ in Section 11.1 found under the assumption that $m = L(I_{2^K+p_k}^{t+K+p})/L(I_{2^K j}^{s+K})$. Arguments for the assumption $m = L(I_{2^K k}^{t+K})/L(I_{2^K+p_j}^{s+K+p})$ are not significantly different. We have

$$d = h\left(\frac{j}{2^s} + \frac{i}{2^{s+K}}\right) = h\left(\frac{j}{2^s} + \frac{i/2^\zeta}{2^{s+K-\zeta}}\right)$$

where $\zeta = z(i)$ and $\lambda(d) = s + K - \zeta$. Since $1 \leq i < 2^K$, we have $K - \zeta > 0$ and $i/2^\zeta$ is an integer with $1 \leq i/2^\zeta < 2^{K-\zeta}$. From this we get $\lambda(d) > s \geq K$ and also

$$\left(\frac{j}{2^s} + \frac{(i/2^\zeta)+1}{2^{s+K-\zeta}}\right) \leq \left(\frac{j}{2^s} + \frac{2^{K-\zeta}}{2^{s+K-\zeta}}\right) = \left(\frac{j+1}{2^s}\right).$$

We have
$$\left(\frac{j}{2^s} + \frac{i/2^\zeta}{2^{s+K-\zeta}}\right) = \left(\frac{j2^{K-\zeta} + i/2^\zeta}{2^{s+K-\zeta}}\right)$$
so that if we let $s' = s + K - \zeta$ and $j' = j2^{K-\zeta} + i/2^\zeta$, then we have $d = x_{j'}^{s'}$, and
$$e = x_{j'+1}^{s'} \leq x_{j+1}^s = b.$$
Since $f$ is affine on $[a,b]$ and $[d,e] \subseteq (a,b)$, we have that $f$ is affine on $[d,e]$.

We have
$$f(d) = h\left(\frac{k}{2^t} + \frac{i/2^\zeta}{2^{s+K+p-\zeta}}\right)$$
and $\lambda(f(d)) = t + K + p - \zeta > t \geq K$. The only point left to argue is that the slope of $f$ to the right of $d$ is not the natural slope of the pair $(d, f(d))$. We note that if we let $t' = t + K + p - \zeta$ and $k' = k2^{K+p-\zeta} + i/2^\zeta$, then $f(d) = x_{k'}^{t'}$.

From the discussion in Section 10.2, we have that the natural slope of a pair of points $(x_\beta^\alpha, x_\delta^\gamma)$ is the ratio $L(I_\delta^\gamma)/L(I_\beta^\alpha)$ as long as $\alpha - \lambda(x_\beta^\alpha) = \gamma - \lambda(x_\delta^\gamma) \geq K$. Note that $\alpha - \lambda(x_\beta^\alpha) = z(\beta)$ and $\gamma - \lambda(x_\delta^\gamma) = z(\delta)$ where $z()$ is the "number of trailing zeros" function defined in Section 11.1.

Let $q$ be a non-negative integer. We have
$$d = h\left(\frac{j}{2^s} + \frac{2^q i}{2^{s+K+q}}\right), \quad \text{and} \quad f(d) = h\left(\frac{k}{2^t} + \frac{2^q i}{2^{t+K+p+q}}\right).$$
Since $d$ is the leftmost point in $[a,b]$ at which $h^{-1}fh$ has a break, we know, for every $r$ with $1 \leq r < 2^q i$, that
$$(11\text{-}4) \qquad \Sigma(x_{\beta(r)}^\alpha) = \Sigma(x_{\delta(r)}^\gamma).$$
From (11-4) and (11-2), we know, for every $r$ with $0 \leq r < 2^q i$, that

(11-5) $\qquad f$ carries $I_{\beta(r)}^\alpha$ affinely onto $I_{\delta(r)}^\gamma$

where $\alpha = s + K + q$, $\gamma = t + K + p + q$, $\beta(r) = 2^{K+q}j + r$ and $\delta(r) = 2^{K+p+q}k + r$ so that
$$x_{\beta(r)}^\alpha = h\left(\frac{j}{2^s} + \frac{r}{2^{s+K+q}}\right), \quad \text{and} \quad x_{\delta(r)}^\gamma = h\left(\frac{k}{2^t} + \frac{r}{2^{t+K+p+q}}\right).$$
We know that (11-4) fails when $r = 2^q i$, so (11-5) must fail for this value of $r$ as well.

We know that $L(I_{\delta(r)}^\gamma)/L(I_{\beta(r)}^\alpha)$ is the natural slope of the pair $(x_{\beta(r)}^\alpha, x_{\delta(r)}^\gamma)$ as long as $z(\beta(r)) = z(\delta(r)) \geq K$. Since $1 \leq i < 2^K$, those values of $r$ with $1 \leq r \leq 2^q i < 2^{K+q}$ have $z(\beta(r)) = z(r)$ and $z(\delta(r)) = z(r)$. We now let $q = K$ and $r = 2^q i = 2^K i$. For this $r$, we have $z(r) \geq K$ so that $L(I_{\delta(r)}^\gamma)/L(I_{\beta(r)}^\alpha)$ is the natural slope of the pair $(x_{\beta(r)}^\alpha, x_{\delta(r)}^\gamma) = (d, f(d))$. But (11-5) fails for this value of $r$ and $f$ is affine on all of $I_{\beta(r)}^\alpha \subseteq [d, e]$. Thus we know that the slope of $f$ to the right of $d$ must be different from the natural slope of the pair $(d, f(d))$. This completes the proof of Theorem 3.

PART V. EXAMPLES

In the next two sections we give examples that illustrate the necessity of the various hypotheses. The first section sets up the machinery by which we build the examples and the second section gives the examples themselves.

## 12. Machinery.

It is very easy to build examples. We essentially give a converse to Section 7. From Lemma 6.4 and the proof of Lemma 6.3 we get the following half of Lemma 6.3.

LEMMA 12.1. *Let $h : S_{n-1} \to S_{n-1}$ be an orientation preserving homeomorphism for which $h(\mathbf{Z}[\frac{1}{n}]) \subseteq \mathbf{Z}[\frac{1}{n}]$, and assume that $h\nu_n h^{-1}$ is in $\overline{T}_{n,n-1}$. Then $hBT_{n,n-1}h^{-1} \subseteq BT_{n,n-1}$.*

The next step is to build homeomorphisms $h$ that conjugate $\nu_n$ into $\overline{T}_{n,n-1}$. However, $\nu_n$ is an expanding map, so that building its conjugate is almost equivalent to building the conjugator. Let $g$ be an element of $\overline{T}_{n,n-1}$ for which there is a Markov partition $P$ that has the following properties.

(1) $P$ is a collection $\{I_0, I_1, \ldots, I_{p-1}\}$ closed intervals of $S_n$ whose union is $S_n$.
(2) The $I_i$ are disjoint except for their endpoints.
(3) Regarding the subscripts modulo $p$, we have that
    (a) each $I_i$ intersects only $I_{i-1}$ and $I_{i+1}$,
    (b) the $I_i$ are arranged in counter-clockwise order around $S_n$ as the subscripts increase, and
    (c) $g$ carries each $I_i$ affinely onto $I_{ni} \cup I_{ni+1} \cup \cdots \cup I_{ni+n-1}$.

We call such a partition of $g$ an affine, $n$-ary Markov partition for $g$. We let $x_i = I_{i-1} \cap I_i$. We have the following.

LEMMA 12.2. *Let $g$ be an element of $\overline{T}_{n,n-1}$ that has an affine, $n$-ary Markov partition $P$. Then there is an orientation preserving homeomorphism $h$ for which $g = h\nu_n h^{-1}$. If we know that the vertices of $P$ are in $\mathbf{Z}[\frac{1}{n}]$ and the number of intervals in $P$ is divisible by $(n-1)$, then $h(\mathbf{Z}[\frac{1}{n}]) \subseteq \mathbf{Z}[\frac{1}{n}]$.*

PROOF: This is standard. Let $p$ be the number of intervals in $P$. We define a Markov partition $Q$ for $\nu_n$ that has $p$ intervals by setting $y_i = (n-1)/pi$ for $o \leq i < p$ and setting $J_i = [y_{i-1}, y_i]$ with subscripts taken modulo $p$. We form derived partitions of $P$ by using $g$ and $Q$ by using $\nu_n$ as discussed in Section 7. It is clear that the endpoints in the derived partitions of $Q$ are dense in $S_{n-1}$.

Let $(a, b)$ be an open interval in $S_{n-1}$. To show that an endpoint of a derived partition of $P$ is in $(a, b)$, we must show that some image $g^k(a, b)$ of $(a, b)$ contains an endpoint of $P$. If $(a, b)$ does not already contain an endpoint of $P$, then $(a, b)$ lies in some $I_i$ of $P$. The length of $(a, b)$ equals $rL(I_i)$ for some real $r$ with $0 < r \leq 1$. If $g(a, b)$ contains no endpoint of $P$, then $g(a, b)$ lies in some $I_j$. Since $g$ is affine on $I_i$, the length of $g(a, b)$ equals $r'L(I_j)$ with $r < r' \leq 1$ and the ratio $r'/r$ depends only on $i$ and $j$. There are only finitely many pairs of intervals from $P$ and so each application of $g$ increases the fraction of an interval in $P$ occupied by the image of $(a, b)$ by an amount that is at least some fixed $\alpha > 1$. Eventually some image of $(a, b)$ under a power of $g$ must contain an endpoint of $P$. Thus the endpoints in the derived partitions of $P$ are dense in $S_{n-1}$.

Now we take endpoints of the derived partitions of $Q$ to the corresponding endpoints of the derived partitions of $P$ and we have built a function $h$ on a dense subset of $S_{n-1}$. Since the endpoints of the derived partitons of $P$ are dense in $S_{n-1}$, the function $h$ extends to a continuous function that has a continuous inverse.

We know that the element $g$ of $\overline{T}_{n,n-1}$ preserves the sets $\mathbf{Z}[\frac{1}{n}]$, $\mathbf{Q}-\mathbf{Z}[\frac{1}{n}]$ and $S_{n,n-1}-\mathbf{Q}$. If we now assume that the endpoints of $P$ are in $\mathbf{Z}[\frac{1}{n}]$, then we have that all endpoints of the derived partitions of $P$ are in $\mathbf{Z}[\frac{1}{n}]$. If the number of intervals in $P$ is divisible by $(n-1)$, then the fixed points of $\nu_n$ (the integer points in $S_{n-1}$) are in the endpoints of $Q$. The elements of $\mathbf{Z}[\frac{1}{n}]$ are precisely the preimages of the integer points in $S_{n-1}$ under iterates of $\nu_n$ and thus must show up in the endpoints of the derived partitions of $Q$. This will give that $h$ takes $\mathbf{Z}[\frac{1}{n}]$ into $\mathbf{Z}[\frac{1}{n}]$.

We now see how to build examples. We will build examnples for $n = 2$. We will build functions $g$ that are conjugate to the doubling map $\nu_2 : S_1 \to S_1$ and that are elements of $\overline{T}_{2,1}$. Since we wish to have the conjugating homeomorphism $h$ have $h(\mathbf{Z}[\frac{1}{2}]) \subseteq \mathbf{Z}[\frac{1}{2}]$, we will build $g$ from an affine, binary Markov partition with $(2-1)2^m = 2^m$ intervals whose endpoints are in $\mathbf{Z}[\frac{1}{2}]$. From Lemma 12.1, we will know that the conjugator $h$ will satisfy $hBT_{2,1}h^{-1} \subseteq BT_{2,1}$. Since $\Delta_2 = \mathbf{Z}[\frac{1}{2}]$, we know from Lemma 2.7 that $BT_{2,1} = T_{2,1}$, so that we will have $hT_{2,1}h^{-1} \subseteq T_{2,1}$. From the proof of half of Lemma 4.3, we will have that a lift $\tilde{h}$ of $h$ to $\mathbf{R}$, will satisfy $\tilde{h}BPL_n(\mathbf{R})\tilde{h}^{-1} \subseteq BPL_n(\mathbf{R})$. From Lemma 2.5, we get that $\tilde{h}\ker(d_n)\tilde{h}^{-1} \subseteq \ker(d_n)$. We note that the lift $\tilde{h}$ of $h$ commutes with the translation $x \mapsto x + (n-1)$ on $\mathbf{R}$. This fact, the fact that $F_n \subseteq \ker(d_n)$ and the definition of $F_n$ imply that $\tilde{h}$ will satisfy $\tilde{h}F_n\tilde{h}^{-1} \subseteq F_n$.

The easiest way to give an affine, binary Markov partition of $S_1$ is to give the sequence of lengths $\{L(I_i)\}$. Since we know that these lengths must add up to 1, we can give these lengths in arbitrary units with the understanding that the final lengths must be normalized to add up to 1. Note that this determines $g$ only up to conjugation by a rotation since we do not say what the values of the endpoints are. We adopt the convention that the endpoint $x_0$ that is the intersection of the first and last intervals is the point 0 in $S_1$. This restricts the set of examples, but not the qualities obtained since conjugation by a rotation does not change any of the properties that we are interested in.

For example, the sequence

(12-1) $$2, 2, 3, 1, 4, 2, 1, 1, 2, 2, 3, 1, 2, 2, 2, 2$$

defines a sequence of 16 intervals lengths. The sum of the lengths given is 32 so that the true length of each interval is $1/32$ times the lengths given. From this it follows that all the endpoints are in $\mathbf{Z}[\frac{1}{2}]$. It is easily checked that the function $g$ defined by this partition is in $\overline{T}_{2,1}$. For example, the first interval $I_0$ is of length 2 (in units of $1/32$) and is carried to the union $I_0 \cup I_1$ of the first two intervals which has length 4. Thus the slope of $g$ on $I_0$ is 2. Similarly

$$\frac{L(I_2 \cup I_3)}{L(I_1)} = \frac{L(I_4 \cup I_5)}{L(I_2)} = \frac{L(I_6 \cup I_7)}{L(I_3)} = 2$$

so that $g$ has slope 2 on $I_1$, $I_2$ and $I_3$, but $L(I_8 \cup I_9)/L(I_4) = 1$ so that $g$ has slope 1 on $I_4$. This gives a break of $-1$ at $x_4 = I_3 \cap I_4$. By summing lengths, we have $x_4 = 8/32 = 1/4$. Further calculations give $g$ breaks of $+1$ at $x_5 = 3/8$ and $x_6 = 7/16$ and $-1$ at $x_8 = 1/2$. The values of $\Sigma()$ are $-2, 0, -1, -1, -2, 0, -2, -1$ at the endpoints of the stable level $x_1, x_3, x_5, x_7, x_9, x_{11}, x_{13}, x_{15}$ respectively.

The fact that the stable level does not have a constant value for $\Sigma()$ implies that $h$ is not PL. This is a much easier observation to make than by all the calculations given above. We leave the next lemma as a very easy exercise for the reader.

LEMMA 12.3. *Let $g \in \overline{T}_{2,1}$ have an affine, binary Markov partition $P$ with $2k$ intervals. The homeomorphism $h$ that conjugates $\nu_2$ to $g$ is PL if and only if for every $i$ with $0 \le i < k$ we have $L(I_{2i}) = L(I_{2i+1})$.*

We can refer to the condition in the lemma as the equal pairs condition. It is easier to state than the condition for the piecewise linearity of $h$ in Section 10.1, but we found it harder to combine with the criterion for the piecewise linearity of $h^{-1}\nu_n h$ in Section 10.2. It turns out that there is a nice inductive proof using the expanding properties of the doubling map on the circle that the criterion of Section 10.1 (for $n = 2$) implies the equal pairs condition. For $n = 2$, this might give an easier argument than that of Section 10.1 that the criterion there implies the piecewise linearity of $h$ (it would eliminate the need for the calculation that verifies (10-1)), but it would hide the rather direct construction of a PL homeomorphism given the criterion. We do not know how much the arguments change when $n = 2$ is not assumed.

## 13. The examples.

EXAMPLE 1: Most of the analysis of Part IV was discovered by studying this example. We use the partition defined by the sequence of lengths in (12-1). This gives an example of a non-PL homeomorphism $h$ that conjugates $\nu_2$ into $\overline{T}_{2,1}$, and, by the remarks in the previous section, conjugates $T_{2,1}$ into $T_{2,1}$ and has lifts that conjugate $F_2$ into $F_2$. In addition, the function $g = h\nu_2 h^{-1}$ has zero break at its fixed point 0. This makes the function $\Sigma()$ well defined. Since

(13-1), $$g^b(x) = \Sigma(x) - \Sigma(g(x))$$

it also makes the non-linearity of $g$ a "coboundary." Thus this example shows that the existence of the function $\Sigma()$ and the equation (13-1) are not sufficient to reach the conclusion of Theorem 3. We note that the right hypotheses to replace (13-1) and the existence of $\Sigma()$ seems to be the assumption that the following stamement fails for only finitely many pairs of points $x$ and $y$ in $S_{n-1}$:

(∗) If $x$ and $y$ have $g^p(x) = g^q(y)$ for some positive integers $p$ and $q$, then $(g^p)^b(x) = (g^q)^b(y)$.

EXAMPLE 2: This example is given by the sequence of 6 lengths

$$2, 2, 1, 1, 1, 1$$

which totals to $8 = 2^3$. Since it satisfies the equal pairs condition, the conjugating homeomorphism $h$ is PL. Since the number of intervals is not a power of 2, there are points in the endpoints of the derived partitions of $Q$ that are not in $\mathbf{Z}[\frac{1}{2}]$. Thus we will have $h(\mathbf{Z}[\frac{1}{2}]) \subseteq \mathbf{Z}[\frac{1}{2}]$, but not $h(\mathbf{Z}[\frac{1}{2}]) = \mathbf{Z}[\frac{1}{2}]$. Thus we have a PL homeomorphism $h$ that conjugates $T_{2,1}$ into $T_{2,1}$ but cannot conjugate $T_{2,1}$ onto $T_{2,1}$. The inverse of $h$ must conjugate $T_{2,1}$ onto a group of PL homeomorphisms of $S_1$ that properly contains $T_{2,1}$.

REMARK: We do not have an example of a non-PL homeomorphism $h$ so that $h$ conjugates $\nu_2$ into $\overline{T}_{2,1}$ and $h^{-1}$ conjugates $\nu_2$ to a PL function. Thus we do not know whether the hypothesis that $h^{-1}\nu_2 h$ be in $\overline{T}_{2,1}$ in Theorem 3 can or cannot be replaced by the hypothesis that $h^{-1}\nu_2 h$ be PL.

EXAMPLE 3: This example is given by the sequence of 16 lengths

$$1, 1, 3, 1, 2, 4, 1, 1, 1, 1, 6, 2, 2, 2, 2, 2$$

which total to $32 = 2^5$. Here the conjugator $h$ fails to have $h(\mathbf{Z}[\frac{1}{2}]) = \mathbf{Z}[\frac{1}{2}]$ in spite of the fact that there are $2^4$ intervals. It is easy to check that the midpoint of $I_5$ (the unique interval of length 4) is a point of period 2 under the action of $g$. This is the point $10/32 = 5/16$ and cannot be the image of any element of $\mathbf{Z}[\frac{1}{2}]$ since the only points of period 2 under the action of $\nu_2$ are the points $1/3$ and $2/3$.

EXAMPLE 4: This example is given by the sequence of 8 lengths

$$1, 3, 4, 2, 1, 3, 1, 1$$

which total to $16 = 2^4$. This example has no well defined function $\Sigma()$ since the break at the fixed point is not zero.

EXAMPLE 5: This example is given by the sequence of 18 lengths

$$1, 1, 3, 1, 2, 1, 3, 1, 2, 2, 1, 3, 2, 1, 1, 3, 2, 2$$

which total to $32 = 2^5$. The points $x_6$ and $x_{12}$ form a cycle of period 2 and $g^b(x_6) = -2$ while $g^b(x_{12}) = +1$. Even though the $g^b(0) = 0$, we cannot define $\Sigma()$ because the sum of the breaks on a future orbit that contains $x_6$ and $x_{12}$ does not converge.